\documentclass[12pt]{article}

\usepackage{amsmath, amssymb, eucal, amscd, color, amsthm, graphicx, mathabx}
\usepackage[all]{xy}

\newcommand{\ca}{\mathcal}

\newcommand{\QQ}{{\mathbb{Q}}}

\newcommand{\ZZ}{{\mathbb{Z}}}
\newcommand\Hom{\operatorname{Hom}}
\newcommand\Ext{\operatorname{Ext}}

\newcommand{\Ker}{\operatorname{Ker}}

\newcommand{\cA}{{\cal A}} \newcommand{\cB}{{\cal B}}
\newcommand{\cD}{{\cal D}}
\newcommand{\cH}{{\cal H}}\newcommand{\cI}{{\cal I}}

\newcommand{\cL}{{\ca L}} \newcommand{\cM}{{\ca M}} 

 \newcommand{\cS}{{\ca S}}\newcommand{\cT}{{\ca T}}
\newcommand{\cU}{{\cal U}}

\newcommand{\al}{\alpha}   
 \newcommand{\eps}{\epsilon}

\newcommand{\scirc}{\circ} 

\newcommand{\ts}{\otimes}

\newcommand{\injto}{\hookrightarrow}
\newcommand{\surjto}{\twoheadrightarrow}
\newcommand{\isoto}{\overset{\sim}{\to}}

\newcommand{\ctop}{\overset{\scirc}}

\newcommand{\mapright}[1]{%
  \smash{\mathop{%
    \hbox to 1cm{\rightarrowfill}}\limits^{#1} } } 
\newcommand{\mapr}{\mapright}
\newcommand{\injmapr}[1]{%
  \smash{\mathop{%
    \hbox to 1cm{\twoheadrightarrowfill}}\limits^{#1} } }

\newcommand{\smapr}[1]{%
  \smash{\mathop{%
    \hbox to 0.5cm{\rightarrowfill}}\limits^{#1} } } 
\newcommand{\mapleft}[1]{%
  \smash{\mathop{%
    \hbox to 1cm{\leftarrowfill}}\limits^{#1} } }
\newcommand{\mapl}{\mapleft}
\newcommand{\maplb}[1]{%
  \smash{\mathop{%
    \hbox to 1cm{\leftarrowfill}}\limits_{#1} } }
\newcommand{\mapdown}[1]{\Big\downarrow
  \llap{$\vcenter{\hbox{$\scriptstyle#1\, $}}$ }}
\newcommand{\mapd}{\mapdown}
\newcommand{\mapdownr}[1]{\Big\downarrow
  \rlap{$\vcenter{\hbox{$\scriptstyle#1\, $}}$ }}
\newcommand{\mapdr}{\mapdownr}
\newcommand{\mapup}[1]{\Big\uparrow
  \llap{$\vcenter{\hbox{$\scriptstyle#1\, $}}$ }}
 \newcommand{\mapu}{\mapup}
\newcommand{\mapupr}[1]{\Big\uparrow
  \rlap{$\vcenter{\hbox{$\scriptstyle#1\, $}}$ }}
 \newcommand{\mapur}{\mapupr}

\newcommand{\ph}{\phantom}

\newcommand{\ti}{\tilde}

\newcommand{\bpf}{\begin{proof}}
 \newcommand{\epf}{\end{proof}}
 \newcommand{\sskp}{\smallskip}
 \newcommand{\bskp}{\bigskip}
 \newcommand{\Gam}{\Gamma}

\newcommand{\simp}{{\rm simp}}

\begin{document}

\newcommand{\TC}{{{}^{\scriptstyle T}\!C}}
\newcommand{\tiS}{\ti{S}}
\newcommand{\chS}{\check{S}}
\newcommand{\cpt}{Cpt}
\newcommand{\chC}{\check{C}}
\newcommand{\tiC}{\tilde{C}}
\newcommand{\fS}{\mathfrak{S}}
\newcommand{\fL}{\mathfrak{L}}

\newcommand{\PropSflabby}{2.1}
\newcommand{\PropBMsoftresolution}{2.2}
\newcommand{\PropgammaZDL}{2.3}

\newcommand{\lowerScalS}{3.1}
\newcommand{\Sectioncsoftsheaf}{3.2}

\newcommand{\Defcoresolution}{3.3}
\newcommand{\Dualcoresolution}{3.4}
\newcommand{\Dualesolution}{3.5}
\newcommand{\doubledualqisom}{3.6}
\newcommand{\cptqisPhiqis}{3.7}
\newcommand{\fScoreslution}{3.8}
\newcommand{\lowerSbiduality}{3.9}

\newcommand{\CheckCptQis}{4.1}
\newcommand{\dualThetaQis}{4.2}
\newcommand{\ThmlfandBM}{4.3}
\newcommand{\xiforS}{4.4}
\newcommand{\Sheafsingularcomparison}{4.5}
\newcommand{\camparisonofisomorphisms}{4.6}

\newcommand{\SectcomparisonBMsupp}{5}
\newcommand{\chSXUtochSZ}{\SectcomparisonBMsupp.1}
\newcommand{\chSXUnatural}{\SectcomparisonBMsupp.2}
\newcommand{\suppcapprod}{\SectcomparisonBMsupp.3}
\newcommand{\comparisonprep}{\SectcomparisonBMsupp.4}
\newcommand{\supportedSheafsingularcomparison}{\SectcomparisonBMsupp.5}

\newcommand{\Sectsimplicial}{6}
\newcommand{\simpchains}{\Sectsimplicial.1}
\newcommand{\suppcapcupprod}{\Sectsimplicial.2}
\newcommand{\xiCDC}{\Sectsimplicial.3}
\newcommand{\xiCS}{\Sectsimplicial.4}
\newcommand{\CScommutes}{\Sectsimplicial.5}
\newcommand{\CScapproduct}{\Sectsimplicial.6}

\newcommand{\SectBMpairs }{7}
\newcommand{\StoSY}{\SectBMpairs.1}
\newcommand{\chSrelative}{7.2}
\newcommand{\ThmlfandBMpair}{7.3}
\newcommand{\ThmlfandBMpaircapprod}{7.4}
\newcommand{\ThmlfandBMpairsimp}{7.5}

\newcommand{\tiSUSUX }{\SectBMpairs.6}
\newcommand{\localisom }{\SectBMpairs.7}
\newcommand{\ThmlfandBMpairlocal}{7.8}
\newcommand{\Thmlocalcapproduct}{7.9}

\title{\bf Borel-Moore homology and cap product operations}
\author{Masaki Hanamura}
\date{}
\maketitle 

\begin{abstract} 
We show that, for a simplicial complex,  the supported cap product operation on 
Borel-Moore homology coincides with the supported  cap product
on simplicial homology.   For this we introduce the supported cap product for 
locally finite singular homology, and compare the cap product on the three homology
theories. 
\end{abstract}

\renewcommand{\thefootnote}{\fnsymbol{footnote}}

In this paper we shall compare the supported cap product operations on the Borel-Moore homology, 
locally finite singular homology, and locally finite simplicial homology; we verify that they agree 
via the canonical identifications of the three homology theories. 

Let $X$ be a locally compact Hausdorff space, and denote by $H_*(X)=H_*(X; \ZZ)$ its Borel-Moore homology 
with $\ZZ$-coefficients. 
The definition (see [BM] and [Br])
 is by means of sheaf theory: take an injective resolution $\cI^*$ of the sheaf $\ZZ$ on $X$,  take its 
``dual" $\cD(\cI^*)$, and then apply the global section functor to get a complex $\Gamma(X, \cD(\cI^*))$; its cohomology 
is the Borel-Moore homology.
For $Z\subset X$ closed, the supported cap product for the Borel-Moore homology 
is a map  
$$H_m(X)\ts H^p_Z(X)\to H_{m-p}(Z)\eqno{(1)}$$
defined via sheaf theory.  In this paper we address the following problems:
\smallskip

(I) \, Suppose $X$ is equipped with a (sufficiently fine) triangulation such that $Z$ is a subcomplex.
Let $H^\simp_*(X)$ denote the locally finite simplicial homology of $X$, defined to be the homology of 
the complex $\tiC_*(X)$ of infinite simplicial chains on $X$.
We have the cap product for  locally finite simplicial homology,
$$H^\simp_m(X)\ts H^p_Z(X)\to H^\simp_{m-p}(Z)\,.\eqno{(2)}$$
It is induced from the chain level map, given
by the cap product formula of ``Alexander-Whitney" type.  
(More precisely, for the chain level map, one needs to replace $Z$ with its closed star -- see \S 6 for details.)
It is well-known that the group $H^\simp_*(X)$ is independent of the choice of 
triangulation. (See [Mu], \S 18 or \S 34 for the homology of finite chains; the same holds for locally finite 
homology.)
Therefore one may ask if this map is independent of the triangulation of $X$ such that $Z$ is a subcomplex. 
\bigskip

(II)\, For a reasonable locally compact Hausdorff space $X$, 
let $H^{lf}_*(X)$ denote the homology of locally finite singular chains of $X$.  
It is known that there is an isomorphism $H_*(X)\cong H_*^{lf}(X)$ (e.g. [Br], Theorem (12.20).)
If $X$ is triangulated, it is well-known that there is a natural isomorphism $H_*^{lf}(X)\cong H_*^\simp(X)$. 
It is thus natural to ask if one can construct a cap product map 
$$H^{lf}_m(X)\ts H^p_Z(X)\to H^{lf}_{m-p}(Z)\,,\eqno{(3)}$$
using only singular chain methods; further one should have:
\smallskip

(i) it is is compatible with (1) via $H_*(X)\cong H_*^{lf}(X)$, and 

(ii) it is compatible with (2) when $X$ is triangulated. 
\smallskip
 
\noindent  
Notice that the construction of cap product (3) together with property (ii) solves the question (I), 
 since cap product (3) is a topological invariant.
\bigskip

We shall answer these questions; it will turn out that the essential problem is the case $Z\neq X$. 
A conclusion of this paper, the topological invariance of the simplicial supported cap 
product (see Corollary to Theorem (\CScapproduct)\,)
apparently has been taken as granted (e.g. see [Bra], section 5), 
 but a proof cannot be found  in the literature.
\bigskip

This paper is organized as follows. 
In \S\S 1 and 2,  we review the definitions of singular homology, cohomology, Borel-Moore homology, and the sheaf theoretic 
supported cap product on Borel-Moore homology.

In \S 3, for a topological space $X$ we 
produce a quasi-isomorphism from the complex of locally finite singular chains $\tiS_*(X)$ to the complex 
$D(\chS^*(X))$, the ``dual'' of the compactly supported cochain complex $\chS^*(X)$. 
If $\cS^*$ is the sheaf of singular cochains, one can take its ``dual" complex $\cD(\cS^*)$, and the complex
$\Gamma(X, \cD(\cS^*))$ computes the Borel-Moore homology. 
It is shown in \S 4 that there is a quasi-isomorphism from $\Gamma(X, \cD(\cS^*))$ to $D(\chS^*(X))$. 
Therefore we have quasi-isomorphisms
$$\tiS_*(X)\to D(\chS^*(X))\gets \Gamma(X, \cD(\cS^*))\,,\eqno{(*)}$$
and on homology an induced isomorphism $H_m^{lf}(X)\to H_m(X)$. 
This gives an explicit, chain-level isomorphism between the 
Borel-Moore homology and the locally finite singular homology.
The presentation $(*)$ itself may not be found in the literature, but one can show 
that it induces the same isomorphism as in [Br], Chap 5, (12.20).

In \S  4, we introduce the singular cap product (3) when $Z=X$;
indeed there is no difficulty there, since the Alexander-Whitney formula provides the correct definition.
The comparison with the sheaf theoretic cap product can be made, in a natural manner, based on the presentation $(*)$.
This section should be contrasted to the next section, which is central in this paper.

In \S 5, we take a closed set $Z$ in general and introduce the cap product 
(3). 
Given $u\in S^p(X, X-Z)$ a cocycle, the Alexander-Whitney formula does not give a map 
a chain map $\cap u: \tiS_*(X)\to\tiS_{*-p}(Z)$. 
This we overcome by resorting to the quasi-isomorphism $\tiS_*(X)\to D(\chS^*(X))$; we attempt to create a chain map 
$u\cup: \chS^{*-p}(Z)\to \chS^*(X)$, which we will indeed do in the derived category. 
After the singular cap product is thus introduced, the comparison with the sheaf theoretic cap product 
again exploits the presentation $(*)$.

In \S 6, the space is assumed to be triangulated.
We first explain the definition of supported cap product in the simplicial homology, 
essentially given on chain level by the Alexander-Whitney formula. 
We show that it is compatible with the singular cap product;
the quasi-isomorphism $\tiS_*(X)\to D(\chS^*(X))$ and its simplicial analogue is 
a necessary component in the proof.

In \S 7, we generalize the results to pairs of spaces. 
If $Y$ is a closed set of the space $X$, then one has the relative 
Borel-Moore homology $H_*(X, Y)=H_*(X, Y;\ZZ)$ defined via sheaf theory. 
Also one has the cap product operation 
$H_m(X, Y)\ts H^p_Z(X)\to H_{m-p}(Z, Z\cap Y)$. 
All of the results of \S\S 4-6 are generalized to this setting, the proofs following the same pattern.
In addition, we discuss the localization isomorphism $H_*(X, Y)\cong H_*(X-Y)$
and its compatibility with cap products.

Throughout the paper, the use of sheaf theory is prevalent. 
Some statements involve no sheaves, but the proofs effectively use them, 
for example in (\lowerSbiduality) and (\chSXUtochSZ). 
\bigskip

The assumption on the space are precisely made as follows:
In \S 1, $X$ is a locally compact Hausdorff space;
in \S 2, we assume also that $X$ is locally contractible and satisfies the second axiom of countability. 
In \S3- \S 5, we further add the assumption that $X$ is of finite cohomological dimension. 
In \S \Sectsimplicial, $X$ is a locally finite, countable simplicial complex of finite dimension.
\bigskip 

We would like to warmly thank Professor T. Suwa for helpful discussions.    
\bigskip

{\bf \S 1.  Singular homology and cohomology }
\bigskip 

For this section we refer the reader to  [Ha], [Mu], and [Sp]. 
Let $X$ be a locally compact Hausdorff space. 
We denote by $S_*(X)$ (resp. $S^*(X)$) the chain complex (resp. cochain complex) of 
singular chains (resp. singular cochains)  on $X$. 
Recall that  $S_p(X)$ is the free abelian group on singular $p$-simplices on $X$, namely continuous 
maps from the topological $p$-simplex $\Delta^p$ to $X$, thus an element of $S_p(X)$ is a finite
$\ZZ$-linear combination of $p$-simplices, 
$$\sum a_\sigma \sigma \,,$$
where $\sigma$ is a singular $p$-simplex and $ a_\sigma\in \ZZ$;
the boundary map $\partial: S_p(X)\to S_{p-1}(X)$ is defined in the usual manner. 
One has $S^*(X):=\Hom (S_*(X), \ZZ)$, and for $u\in S^p(X)$ and $x\in S_{p+1}(X)$, 
$(du)(x)=(-1)^{p+1} u(\partial x)$. This sign convention differs from the usual one.
The homology of $S_*(X)$ is denoted $H_*^c(X)$, and the cohomology of $S^*(X)$ is denoted 
$H^*(X)$. 

We also recall that, for a subset $A\subset X$, there correspond the relative singular chain complex $S_*(X, A)$ and
the relative singular cochain complex $S^*(X, A)$. 

Let $\tiS_*(X)$ be the chain complex of locally finite chains on $X$; an element of 
$\tiS_p(X)$ is  a possibly infinite sum $\sum a_\sigma \sigma$ with $p$-simplices $\sigma$,
which is locally finite.
One defines $H^{lf}_* (X)$ to be the homology of $\tiS_*(X)$. 
The functor $X\mapsto \tiS_*(X)$, thus also $H^{lf}_* (X)$,  is covariantly functorial for proper maps. 
(Locally finite homology appears as homology of the second kind in Cartan Seminar, 1948-49.)

We define a subcomplex $\chS^*(X)$ of $S^*(X)$ by 
$$\chS^p(X)=\varinjlim S^p(X, X-K)\,.$$
where $K$ varies over the compact subsets of $X$. 
We denote the cohomology of $\chS^*(X)$ by $H^*_c(X)$. 
The functors $X\mapsto \chS^*(X)$ and $H_{c}^*(X)$ are cotravariantly functorial for proper maps.

There is the cup product map 
$\cup :S^*(X)\ts S^*(X) \to S^*(X)$ given for $u\in S^p(X)$, $v\in S^q(X)$ by
$$(u\cup v)([v_0, \cdots, v_{p+q}])=(-1)^{pq}u([v_0, \cdots, v_p])v([v_p, \cdots, v_{p+q}])\,;$$
notice again the difference from the usual convention. 
One has $d(u\cup v)=(du)\cup v+(-1)^p u\cup dv$.

The cap product map 
$\cap : S_*(X)\ts S^*(X)\to S_*(X)$ is given as follows.
For $u\in S^p(X)$ and $\sigma=[v_0, \cdots, v_m]$  a singular $m$-simplex, let 
$$\sigma \cap u =  u(\sigma') \sigma''$$  
where $\sigma'=[v_0, \cdots, v_p]$ and  $\sigma''=[v_p, \cdots, v_m]$.
We have the identity
$$\partial(\al\cap u)=  (-1)^p(\partial \al)\cap u +\al\cap (du)\,.$$
For $\al\in S_m(X)$ and $u\in S^p(X)$, $v\in S^{m-p}(X)$, 
one verifies
$$v(\al\cap u)=(-1)^{p(m-p)} (u\cup v)(\al)\,.$$
\bigskip

{\bf \S 2.  Borel-Moore homology and sheaf theoretic cap product}
\bigskip 

The references for the Borel-Moore homology include [BM], [Br] and [I]; we will use [Br] as our 
main reference. 
In this section, we assume:

\begin{quote} (*) \quad $X$ is a locally compact Hausdorff topological space, which satisfies the second axiom of countability, and 
which is locally contractible.
\end{quote}

\noindent For example a locally finite, countable CW complex satisfies this condition.
We note and use the following facts:
\smallskip 

 $\bullet$\qquad A space $X$  satisfying (*) is paracompact. 
 
 $\bullet$\qquad An open set of $X$ also satisfies  (*).  A closed subset of $X$ satisfies all the conditions of (*) except
 the local contractibility. 
\bigskip 
 
The assignment $U\mapsto S^p(U)$, where $S^p(U)$ is the group of singular $p$-cochains on $U$, 
gives a presheaf on $X$, and let $\cS^p$ be the associated sheaf. Thus we have the differential sheaf 
$\cS^*$, called the singular cochain sheaf on $X$ ([Br- Chap. I, \S 7]).

Since $S^p(U)$ is conjunctive ([Br-I, p. 26]), applying [Br-I, (6.2)], 
one has an exact sequence 
$$0\to S^p_0(U) \to S^p(U) \mapr{\theta} \Gamma(U, \cS^p)\to 0\eqno{}$$
where $\theta$ is the canonical map, and 
$S^p_0(U)$ is the kernel of $\theta$; $S^p_0(U)$ is the complex of singular cochains 
on $U$ that are locally zero. The subcomplex $S^*_0(U)$ is acyclic, so $\theta$ is a quasi-isomorphism. 
\bigskip 

\noindent (2.1) {\bf Proposition.}\quad{\it $\cS^*$ is a resolution of $\ZZ$ on $X$
by flabby sheaves.
}\smallskip 

\begin{proof}  
The flabbiness follows from the surjectivity of $\theta$ and the surjectivity 
of the restriction maps $S^p(X)\to S^p(U)$. 
The exactness of $0\to \ZZ\to \cS^*$ follows from [Br-II-(1.2)], since 
$X$ is assumed locally contractible
(indeed the weaker condition $HLC_\ZZ^\infty$ will suffice.)
\end{proof}

By this fact, the sheaf cohomology of $\ZZ$, $H^p(X; \ZZ)$, is identified with 
$H^p\Gamma(X, \cS^*)$. 
Since $\theta: S^*(X)\to \Gamma(X, \cS^*)$ is a quasi-isomorphism, we also have 
$H^p(S^*(X))\cong H^p\Gamma(X, \cS^*)$.  We will write $H^p(X)$ for $H^p(X; \ZZ)$. 

We also note that for $Z$ a closed subset of $X$, we have $H^p(S^*(X, X-Z))\cong 
H^p_Z(X; \ZZ)$. 
Indeed, there are quasi-isomorphisms $S^*(X)\to \Gamma(X, \cS^*)$ and $S^*(X-Z)\to \Gamma(X-Z, \cS^*)$
which make the right square in the following diagram commute
$$
\begin{array}{ccccccc}
 0\to&   S^*(X, X-Z) &\mapr{}    &S^*(X)     &\mapr{} &S^*(X-Z) &\to 0 \\
&\mapd{\theta}&    & \mapd{}   & &\mapd{}  &\\
0\to &\Gamma_Z(X, \cS^*)&\mapr{}   &\Gamma(X, \cS^*)   &\mapr{} &\Gamma(X-Z, \cS^*) &\ph{\,.} \to 0\,.
\end{array}
$$
Since the rows of the diagram are exact, 
there is a unique map of complexes $\theta: S^*(X, X-Z)\to \Gamma_Z(X, \cS^*)$ so that the 
left square commutes, and it is also a quasi-isomorphism.
\bigskip

For a complex $K^\bullet$ of abelian groups, its dual $D(K^\bullet)$ is defined by
$$D(K^\bullet)=\Hom(K^\bullet, I^\bullet)\,, $$
where $I^\bullet$ is the complex $[\QQ \to \QQ/\ZZ]$ concentrated in degrees 0 and 1. 
For $f\in \Hom(K^\bullet, I^\bullet)$ and $x\in K^\bullet$,  $df$ is defined by the formula
$$(df)(x)=(-1)^{|f|+1} f(dx)+ d(f(x))\,$$
(where $|f|$ denotes the degree of $f$). 
Note that $I^\bullet$ is an injective resolution of $\ZZ$, and $D(K^\bullet)=R\Hom(K^\bullet, \ZZ)$ in the derived category of abelian groups.
The functor $D$ is exact, and takes quasi-isomorphisms to quasi-isomorphisms.
If $K^\bullet$ is a complex of free $\ZZ$-modules, then the natural map 
$\Hom(K^\bullet, \ZZ) \to D(K^\bullet)$ is a quasi-isomorphism. 

If $\cL^*$ is a complex of $c$-soft sheaves on $X$, its dual $\cD(\cL^*)$ is the complex of flabby sheaves given
by 
$$U\mapsto \Hom(\Gamma_c(U, \cL^*), I^\bullet)=D\Gamma_c(U, \cL^*)\,,$$
see [Br-V, \S 2].

Recall from [Br-V, \S 3] that the Borel-Moore homology of $X$ (assumed only locally compact 
Hausdorff) is defined as
$$H_p(X; \ZZ)=H_p(X)=H_p\Gamma (X, \cD(\cI^*))\,,$$
in which $\cI^*$ is the canonical injective resolution of $\ZZ$, and $\cD(\cI^*)$ 
is the dual of $\cI^*$. One should distinguish it from homology with compact support $H^c_p(X; \ZZ)$. 
We have (see [Br, 293, (10)]):
\bigskip

\noindent (2.2) {\bf Proposition.}\quad{\it  If $\cL^*$ is a $c$-soft resolution of $\ZZ$, then one has
$H_p(X; \ZZ)=H_p\Gamma (X, \cD(\cL^*))$. 
}\bigskip 

Now for our $X$, since $\cS^*$ is a flabby resolution of $\ZZ$
(in particular a $c$-soft resolution of $\ZZ$), we have
$$H_p(X; \ZZ)=H_p\Gamma (X, \cD(\cS^*))\,, $$
namely the Borel-Moore homology can be calculated using $\cS^*$.

One has the product map $\cS^*\ts \cS^*\to \cS^*$, that is induced from the 
cup product of singular cochains, and it is compatible with augmentations. 
Thus according to [Br-V-(10.3)]) one can use it to produce a
map of differential sheaves (cap product at the sheaf level) 
$$\cap: \cD(\cS^*)\ts \cS^* \to \cD(\cS^*)\,,$$
where $f\cap s$ for $f$ of degree $m$ and $s$ of degree $p$ is defined by 
$$\langle  f\cap s, t\rangle =(-1)^{mp}\langle f, s\cup t\rangle \,.$$
(The value of a functional $f$ at $x$ will be written $\langle f, x\rangle $. )
Then one has 
$$d(f\cap s)=(-1)^p df\cap s+ f\cap ds\,.$$

For $Z$  a locally contractible closed subset of $X$, one has the induced map on sections
$$\cap: \Gamma(X, \cD(\cS^*))\ts \Gamma_Z(X, \cS^*) \to \Gamma_Z(X, \cD(\cS^*))=\Gamma(Z, \cD(\cS^*|_Z))\,.$$
The last equality holds by the next proposition.
Since $\cS^*|_Z$ is a $c$-soft resolution of $\ZZ$ on $Z$, the homology of the complex $\Gamma(Z, \cD(\cS^*|_Z))$ 
is identified with the Borel-Moore homology $H_*(Z)$. 
Therefore, passing to homology we obtain a map 
$$\cap : H_m(X)\ts H^p_Z(X)\to H_{m-p}(Z)\,;$$
this is the {\it sheaf theoretic} supported cap product. 
\bigskip 

\noindent (2.3) {\bf Proposition.}\quad{\it
If $\cL$ a $c$-soft sheaf on $X$ and $Z$ a closed set, then 
$$\Gamma_Z (X, \cD(\cL))=\Gamma (Z, \cD(\cL |_Z))\,$$
}\smallskip 

\begin{proof} This is a special case of [Br-V-(5.5)], but
we give a direct, simpler proof. 

One has an exact sequence (with $U=X-Z$)
$$0\to \Gamma_c(U, \cL) \to  \Gamma_c(X, \cL)\to \Gamma_c(Z, \cL|_Z)\to 0\,,$$
the last surjection being a consequence of $c$-softness.  
Taking dual, one has an exact sequence
$$0\to \Gamma(Z, \cD(\cL|_Z))\to \Gamma(X, \cD(\cL)) \to \Gamma(U, \cD(\cL)) \to 0\,.$$
The kernel of the second map equals  $\Gamma_Z(X, \cD(\cL))$. 
\end{proof}
\bigskip

{\bf \S 3. The sheaf $\cS_*$ and the complex of locally finite singular chains.}
\bigskip 

For the rest of the paper (indeed starting with (\lowerSbiduality)),  we add another condition to (*) and assume that

\begin{quote}
(**) \quad $X$ is a locally compact Hausdorff topological space which satisfies the second axiom of countability,
which is locally contractible, and $\dim_\ZZ X<\infty$. 
(See  [Br, II-\S 16] for the notion of cohomological dimension.)
\end{quote}

\noindent  
For example a locally finite countable CW complex of 
finite dimension satisfies these conditions. 
If $X$ satisfies (**), then an open set of $X$ also satisfies  (**), and 
a locally contractible closed set of $X$ also satisfies  (**)
(see [Br, II-(16.8)] for the notion of dimension of spaces). 
\bigskip 

We define a map of complexes $\xi: \tiS_*(U) \to \Hom(\chS^*(U), \ZZ)\subset 
D\chS^*(U)$ as follows.
Let  $\al\in \tiS_m(X)$. 
For  $u\in \chS^m(X)$, let $K\in \cpt(X)$ such that
$u\in S^m(X, X-K)$, write $\al=\al'+\al''$ with $\al'\in S_m(X)$, $\al''\in \tiS_m(X-K)$, 
and define $\xi(\al)\in \Hom(\chS^m(X), \ZZ)$ by
$$\langle \xi(\al), u\rangle =(-1)^{m}\langle u, \al'\rangle \,.$$
This is well-defined independent of the choice of $K$ and the decomposition of $\al$;
one also verifies that it gives a map of complexes.

Toward the end of this section we will prove that the map $\xi$ is a quasi-isomorphism; 
it is all we need in later sections. 
The reader may opt to grant it and proceed to \S 4. 
\bigskip

{\bf The complex $\cS_{X, *}$.}\quad 
We recall from [Br] the definition and properties of the complexes of sheaves $\underline{\Delta}_{X, *}$
and $\cS_{X, *}$. 
\smallskip 

1. Consider the presheaf $U\mapsto S_p(X, X-U)$ on $X$, and let $\underline{\Delta}_{X, p}$ be the associated sheaf. 
For any paracompactifying family  of supports $\Phi$, 
$\Gamma_\Phi(X, \underline{\Delta}_{X, p})$ coincides with the group of locally finite singular chains with support 
in $\Phi$;
 indeed by [Br, p.31, Exercise 12] $\underline{\Delta}_{X, p}$ is a monopresheaf and conjunctive for coverings of $X$, 
 so that [Br-I-(6.2)] applies.  Since $X$ is paracompact, we may take $\Phi=cld$ (the family of all closed sets of $X$) and since $X$ is locally compact we may
 take $\Phi=c$. 
\smallskip 

2. Let $\fS_p(X)=\varinjlim A_n$, where $A_n=S_p(X)$ for each $n\ge 1$, and the map $A_n \to A_{n+1}$ is the subdivision map.
There is a natural injection $S_*(X)\to \fS_*(X)$, which is a quasi-isomorphism.
For $A\subset X$, one defines $\fS_*(X, A)=\fS_*(X)/\fS_*(A)$. 
The assignment $U\mapsto \fS_p(U)$, where $U$ is an open set of $X$, is a flabby cosheaf on $X$(see [Br-V, \S 1]). 

In general, if  $\cL$ is a $c$-soft sheaf on $X$, then we denote by $\Gamma_c\{\cL\}$ the flabby cosheaf 
$U\mapsto \Gamma_c(U, \cL)$ (cf. [Br-V-(1.6)]). 
Conversely it is known that a flabby cosheaf $\frak{L}$ 
is of the form $\Gamma_c\{\cL\}$ for a unique $c$-soft sheaf $\cL$, [Br-V-(1.8)]

To the flabby cosheaf $\fS_p(U)$ there corresponds a $c$-soft sheaf $\cS_{X, p}$, 
thus we have $\fS_p(U)=\Gamma_c(U, \cS_{X, p})$. 
From the proof of this correspondence, $\cS_{X, p}$ is 
the sheaf associated with the presheaf  $U\mapsto \mathfrak{S}_p(X, X-U)$.
\smallskip 

3.   One has a natural injection 
$\underline{\Delta}_{X, p}\to \cS_{X, p}$, 
induced from the injection $S_p(X, X-U)\injto \fS_p(X, X-U)$. 
For any paracompactifying family of supports $\Phi$, the induced map 
$$\Gamma_\Phi(X, \underline{\Delta}_{X, *}) \injto \Gamma_\Phi(X,\cS_{X, *})$$
is a quasi-isomorphism (cf. [Br-V-(1.19)]).
In particular, we have quasi-isomorphisms
$$\ti{S}_*(X)=\Gamma(X, \underline{\Delta}_{X, *}) \injto \Gamma(X,\cS_{X, *})$$
and 
$$S_*(X)=\Gamma_c(X, \underline{\Delta}_{X, *}) \injto \Gamma_c(X,\cS_{X, *})\,.$$
\smallskip 

4.  The formation of the sheaf $\cS_{X, p}$ is compatible with restriction to open sets, in the following sense. 
Let $\cS_{U, p}$ be the sheaf on $U$ associated with the  presheaf $V\mapsto \fS_p(U, U-V)$.
The natural injection $\fS_p(U, U-V)\injto \fS_p(X, X-V)$,  for $V\subset U$ open, induces a map 
of sheaves $\cS_{U, p} \to  \cS_{X, p}|_U$.  Since both sheaves are  $c$-soft,  and 
the induced map on each open set
$$\Gamma_c(V, \cS_{U, p}) \to \Gamma_c(V, \cS_{X, p}|_U) $$
can be  identified with the identity map on $\mathfrak{S}_p(V)$,
 it follows that the map 
 $\cS_{U, p} \to  \cS_{X, p}|_U$ is an isomorphism.   
 Because of this fact, we will usually write $\cS_p$ for  $\cS_{X, p}$ or $\cS_{U, p}$. 
Composing the quasi-isomorphism $ \ti{S}_*(U) \injto \Gamma(U,\cS_{U, *})$ with the isomorphism $ \Gamma(U,\cS_{U, *})\cong 
 \Gamma(U,\cS_{X, *})$, we obtain a quasi-isomorphism $\ti{S}_*(U) \injto \Gamma(U,\cS_{X, *})$. 
 Similarly we have a quasi-isomorphism ${S}_*(U) \injto \Gamma_c(U,\cS_{X, *})$. 
 \bigskip 

{\bf Dual complexes of $\ti S_*(X)$.}\quad   We introduce duals and double duals of the complex $\ti S_*(X)$;
instead of naive duals, we pass to a direct family of complexes, and then proceed to take duals. 

For an open set $U$ of $X$, and a compact $K\subset U$, one has a complex 
$S_*(U, U-K)$.  For $K\subset K'$, there is an induced surjective map 
$S_*(U, U-K')\to S_*(U, U-K)$, so one has an inverse system of complexes $\{S_*(U, U-K)\}$, 
where $K$ varies over $Cpt(U)$, the directed set of compact subsets of $U$. 
One has $\tiS_*(U)=\varprojlim_{K\in Cpt(U)} S_*(U, U-K)$, as is easily verified.
Taking the dual,  we have a direct system of complexes $\{DS_*(U, U-K)\}$, and we let 
$$\widecheck{DS}_*(U):=\varinjlim_{K\in Cpt(U)} DS_*(U, U-K)\,. $$
The natural inclusion $S^*(U, U-K)\injto D(S_*(U, U-K))$ induces an injection
$ \check{S}^*(U) \injto \widecheck{DS}_*(U)$. 
Taking another dual, one has  an inverse system  $\{DD S_*(U, U-K)\}$, and 
$$\widetilde{DDS}_*(U):=\varprojlim_{K\in Cpt(U)} DDS_*(U, U-K)\,.$$
\bigskip 

\noindent (3.1) {\bf Proposition.}\quad{\it  Let $\cS_*=\cS_{X, *}$ and let $U$ be an open set of $X$. 
We have maps of complexes, each of them being a quasi-isomorphism,
$$\begin{array}{cccccc}
    & &\ti{S}_*(U)    &\overset{\Theta}\injto &\Gamma(U, \cS_{*})\,,  &\ph{aaaaaaa} (1)\ph{_c}\\
    & &S_*(U)        &\overset{\Theta}\injto &\Gamma_c(U, \cS_{*}) \,,    &\ph{aaaaaaa} (1)_c\\
 S^*(U)   &\injto{}  & DS_*(U)   &\overset{\Theta'}\twoheadleftarrow  &\Gamma(U, \cD(\cS_*))\,, &\ph{aaaaaaa} (2)\ph{_c}\\
  \check{S}^*(U)  &\injto{}  &\widecheck{DS}_*(U)    &\overset{\Theta'}\twoheadleftarrow   &\Gamma_c(U, \cD(\cS_*))\,, &\ph{aaaaaaa} (2)_c\\
 D( \chS^*(U) )  &\twoheadleftarrow{}  &\widetilde{DDS}_*(U)    &\overset{\Theta''}\injto  &\ph{\,} \Gamma(U, \cD\cD(\cS_*))  \,.&\ph{aaaaaaa} (3)\ph{_c}
   \end{array}
$$ 
}\smallskip 

{\bf Remark.}\, The statement is made for any open set $U$ of $X$ to clarify the relevance of sheaf theory; one could have stated it only for
$U=X$ without losing generality.
\smallskip 

\begin{proof}
As already noted, we have a map
$$\Theta:\ti{S}_*(U)\injto \Gamma(U, \cS_{U, *})=\Gamma(U, \cS_{X, *})$$
which is a quasi-isomorphism, verifying (1).
We argue similarly for $(1)_c$. 

From $(1)_c$, we have induced maps
$$S^*(U)\hookrightarrow DS_*(U) \overset{\Theta'}\twoheadleftarrow  D\Gamma_c(U, \cS_*)=\Gamma(U, \cD(\cS_*))\,$$
with both arrows quasi-isomorphisms (note that $S_*(U)$ is free, so the first map is a quasi-isomorphism); this is (2). 

For $K\in \cpt(U)$,  the quasi-isomorphism
$$\Theta: S_*(U, U-K)\hookrightarrow \Gamma_c(U, \cS_*)/\Gamma_c(U-K, \cS_*)$$
induces quasi-isomorphisms
$$S^*(U, U-K)\hookrightarrow DS_*(U, U-K) \overset{\Theta'}\twoheadleftarrow  D(\Gamma_c(U, \cS_*)/\Gamma_c(U-K, \cS_*))
=\Gamma_K(U, \cD(\cS_*))\, $$
and by taking direct limit, we obtain quasi-isomorphisms
$$\check{S}^*(U)\hookrightarrow \widecheck{DS}_*(U)\overset{\Theta'}\twoheadleftarrow  \Gamma_c(U, \cD(\cS_*))\,,$$
giving $(2)_c$.

(3) follows from $(2)_c$ by taking dual, since  $\Gamma(U, \cD\cD(\cS_*))=D\Gamma_c(U, \cD(\cS_*))$ and 
$D(\widecheck{DS}_*(U))=\widetilde{DDS}_*(U)$. 
\end{proof}

To relate the map (1) to the diagram (3)  in Proposition (\lowerSbiduality), we introduce a few maps related to duality.
\smallskip 

1.  For a complex of abelian groups $K$, let 
$\frak{d}: K\to DD(K)$ be the map of complexes which sends 
$x\in K$ to the element $\hat{x} \in DD(K)$ given by
$$\langle \hat{x}, f\rangle =(-1)^{|x|\cdot |f|} \langle f, x\rangle , \quad f\in D(K)\,.$$
Here $|x|$, for example,  denotes the degree of $x$. 

In particular, applying this to the complex $S_*(U, U-K)$
one has the map $\frak{d}: S_*(U, U-K)\to DDS_*(U, U-K)$. 
Passing to the inverse limit
for $K\in \cpt(U)$, we obtain a map 
$$\frak{d}: \tiS_*(U)\to \widetilde{DDS}_*(U)$$
which takes $\al=(\al_K)$ to $\frak{d}(\al)=(\widehat{\al_K})$.
\smallskip 

2. For $\cL_*$ a bounded below complex of $c$-soft sheaves, 
one has a canonical map 
$$\frak{d}:\cL_*\to \cD\cD(\cL_*)$$
 defined as follows. 
Let  $x\in \Gamma(U, \cL_p)$, and we shall define its image $\hat{x}\in \Gamma(U, \cD\cD(\cL_p))
=D\Gamma_c(U, \cD(\cL_p))$. 
Take an element $f\in \Gamma_c(U, \cD(\cL_p))$, 
which may be viewed as an element of a larger group
$\Gamma(U, \cD(\cL_p))=D\Gamma_c(U, \cL_p)$. 
Let $|f|=K$, 
choose $x_1\in \Gamma_c(U, \cL_p)$ such that $x|_K=
x_1|_K\in\Gamma(K, \cL_p)$, and set 
$$\langle \hat{x}, f\rangle =(-1)^{|f|\cdot |x|} \langle f, x_1\rangle \in I^\bullet \,.$$
This is well-defined independent of the choice of $x_1$. 

{\bf Remark.}\quad The map $\frak{d}$ does not appear as such in [Br], but it induces the map in [Br, p.285, (1.13)]
upon taking $\Gamma_c\{-\}$.
\smallskip 

We leave it as an exercise to verify that the composition of the maps $\frak{d}: \tiS_*(U)\to \widetilde{DDS}_*(U)$
and $\widetilde{DDS}_*(U) \to D( \chS^*(U) ) $ coincides with $\xi$. 

These give us a diagram 
$$\begin{array}{ccccc}
&& \ti{S}_*(U) &\mapr{\Theta}  & \Gamma(U, \cS_*) \\
&\raise1ex\hbox{$\scriptstyle{\xi}$}\!\!\swarrow &\mapd{\frak{d}} &     &\mapdr{\frak{d}} \\
D(\chS^*(U))&\mapl{}&\widetilde{DDS}_*(U)   &\mapr{\Theta''}  &\ph{\,.}\Gamma(U, \cD\cD(\cS_*))\,
 \end{array}
$$  
and the square on the right commutes as we now show.
\bigskip

We need an alternative description of the map $\Theta:\ti{S}_*(U) \injto \Gamma(U,\cS_{X, *})$.
In general, for a $c$-soft sheaf $\cL$ on $X$, and an open set $U$, consider the inverse system 
of abelian groups $\Gamma_c(U, \cL)/\Gamma_c(U-K, \cL)$ indexed by $K\in Cpt(U)$. 
We claim that there is a natural map to this inverse system
$$\Gamma(U, \cL) \to \{\Gamma_c(U, \cL)/\Gamma_c(U-K, \cL)\}_K\,.$$
Indeed, for each $s\in \Gamma(U, \cL)$, there is a decomposition $s=s'+s''$ where 
$s'\in \Gamma_c(U, \cL)$ and $s''\in \Gamma(U, \cL)$ with $|s''|\subset U-K$. 
(The restriction map $\Gamma_c(U, \cL)\to \Gamma(K, \cL|_K)$ is surjective since $\cL|_U$ is $c$-soft, 
so one can take $s'\in \Gamma_c(U, \cL)$ such that $s'|_K=s|_K$.)
The desired map $\Gamma(U, \cL) \to \Gamma_c(U, \cL)/\Gamma_c(U-K, \cL)$ is given by $s\mapsto \bar{s'}$. 
\bigskip 

\noindent (3.2) {\bf Proposition.}\quad{\it 
The above map induces 
 an isomorphism to the inverse limit
 $$\Gamma(U, \cL) \isoto \varprojlim ( \Gamma_c(U, \cL)/\Gamma_c(U-K, \cL))\,.$$
 }\bigskip 
 \begin{proof} The injectivity is obvious. 
 For surjectivity, let $s_K\in \Gamma_c(U, \cL)/\Gamma_c(U-K, \cL)$ be a family of elements, indexed by 
 $K$, which is coherent:  $K'\supset K$ implies $s_{K'}-s_K \in \Gamma_c(U-K, \cL)$. 
 Then the family of elements, indexed by $K$, 
 $$s_K|\mathring{K}\in \Gamma(\mathring{K}, \cL)\,,$$
 satisfies the patching condition, so it gives a global section restricting to each $s_K$.
 \end{proof}

Going back to our situation, the map of complexes $S_*(U) \to \Gamma_c(U,\cS_{X, *})$ induces a map of inverse systems
of complexes
$S_*(U, U-K) \to \Gamma_c(U,\cS_{X, *})/\Gamma_c(U-K,\cS_{X, *})$, and taking inverse limit we obtain 
a map $\ti{S}_*(U) \injto \Gamma(U,\cS_{X, *})$.  The reader may verify that this coincides with the map $\Theta$
given before. 
\bigskip 

For $K\in \cpt(U)$ we  have a commutative diagram 
$$\begin{array}{ccc}
 S_*(U, U-K)   &\mapr{\Theta}  &\Gamma_c(U, \cS_*)/ \Gamma_c(U-K, \cS_*) \\
\mapd{\frak{d}} &     &\mapdr{\frak{d}} \\
 DD( S_*(U, U-K) \,)&\mapr{\Theta''}  &\ph{\,.} DD(\Gamma_c(U, \cS_*)/ \Gamma_c(U-K, \cS_*))\,.
 \end{array}
$$  
Passing to the inverse limit, and using the second description of $\Theta$, we obtain 
the commutativity of the right square in the diagram in question.
\bigskip

In (\lowerSbiduality), under a further assumption on cohomological finite dimensionality on $X$, 
it will be proven that the map $\frak{d}$ on the right is a quasi-isomorphism, and hence all the maps in the
diagram are quasi-isomorphisms.
The proof is based on [Br, Chap. V], which takes the rest of the section.
The reader may choose to skip it, granting (\lowerSbiduality) as a fact. 

We need theorems (\doubledualqisom)-(\fScoreslution), essentially from [Br]. 
The corresponding statements in [Br] are somewhat different, due partly to the emphasis on
cosheaves (rather than $c$-soft sheaves) and to the generality on assumptions.  
We shall therefore recall some notions and results from [Br] (in terms of $c$-soft sheaves, when possible),
and give comments on the statement and the proof of (\doubledualqisom).
 
For the definitions of precosheaves, cosheaves, and local isomorphism of precosheaves, see [Br-V, \S 1 and \S 12].
We begin here with the definition of resolutions and coresolutions of the constant sheaf $\ZZ$ on $X$. 
\bigskip 

\noindent (3.3) {\bf Definition.}([Br-V-(12.5)])\quad Let $\cL^*$ be a bounded below complex of sheaves on $X$
and $\eps: \ZZ \to H^0(\cL^*)$
be a homomorphism of presheaves, where $H^0(\cL^*)$ is the presheaf $U\mapsto H^0(\cL^*(U))$. We say
that $(\cL^*, \eps)$ is a {\it quasi-resolution} of $\ZZ$ if  $\cH^p(\cL^*)=0$ for $p\neq 0$ and $\eps$ induces an
isomorphism of sheaves $\eps: \ZZ \to \cH^0(\cL^*)$. 
Here $\cH^p(\cL^*)$ denotes the sheafication of the presheaf $H^p(\cL^*)$. 

Let $\fL_*$ be bounded below complex of cosheaves and $\eta: H_0(\fL_*) \to \ZZ$ be a homomorphism of
precosheaves, where $H_0(\fL_*)$ is the precosheaf $U\mapsto H_0(\fL_*(U))$. 
Then $(\fL_*, \eta)$ is said to be a {\it quasi-coresolution} of $\ZZ$ if $H_p(\fL_*)=0$ for $p<0$, 
$H_p(\fL_*)=0$ is locally zero for $p>0$, and 
$\eta$ is a local isomorphism of precosheaves. 
\bigskip 

The main example of a quasi-coresolution arises from the complex $\cS_*$, see (\fScoreslution). 

Indeed in [Br] more general notions of quasi-$n$-resolution of $\ZZ$ and quasi-$n$-coresolution of $\ZZ$ are
introduced; we have restricted ourselves to the case $n=\infty$. 
We cite below some results from [Br-V]; these theorems are formulated for quasi-$n$-resolutions and quasi-$n$-coresolutions, 
but we again take $n=\infty$.   

 Recall the notation $\Gamma_c\{\cL\}$ for the cosheaf associated with a $c$-soft sheaf $\cL$. 
According to the next two facts, the dual of a coresolution  is a resolution, and the dual of a resolution  is a 
coresolution. 
\bigskip

\noindent (3.4) {\bf Theorem.}([Br-V-(12.7)])\quad{\it 
Let $\cL_*$ be a bounded below complex of $c$-soft sheaves such that $\Gamma_c\{\cL_*\}$
is a quasi-coresolution of $\ZZ$.  Then $X$ is $clc_\ZZ^\infty$ and $\cD(\cL_*)$ is a quasi-resolution of $\ZZ$
by flabby sheaves. 
}\bigskip 

\noindent (3.5) {\bf Theorem.}([Br-V-(12.9)])\quad{\it
Suppose that $X$ is $clc_\ZZ^\infty$. 
If $\cL^*$ is a $c$-soft quasi-resolution of $\ZZ$, then the complex of cosheaves 
$\Gamma_c\{\cD(\cL^*)\}$ is a quasi-coresolution of $\ZZ$. 
}\bigskip 

For the notion of cohomological local connectivity $clc_\ZZ^\infty$, see [Br-II, \S 17] (although its
knowledge is not needed for what follows). 
The proof of the following uses (\Dualcoresolution) and (\Dualesolution). 
\bigskip 

\noindent (3.6) {\bf Theorem.}\quad{\it  
If $\cL_*$ is a bounded below complex of $c$-soft sheaves such that $\Gamma_c\{\cL_*\}$
is a quasi-coresolution of $\ZZ$, then the canonical map 
$$\frak{d}: \cL_*\to \cD\cD(\cL_*)$$
induces a quasi-isomorphism 
$$\Gamma_c(U, \cL_*) \to \Gamma_c(U, \cD\cD(\cL_*))$$
for each open set $U$ of $X$. 
}\smallskip 

\begin{proof}
This is essentially due to  [Br-V-(12.11)] particularly its proof.  By  (\Dualcoresolution), $X$ is $clc_\ZZ^\infty$ and
$\cD(\cL_*)$ is a quasi-resolution of $\ZZ$ by flabby sheaves.
Then by (\Dualesolution),  $\Gamma_c\{\cD\cD(\cL^*)\}$ is a quasi-coresolution of $\ZZ$. 
The map $\frak{d}:\cL_*\to \cD\cD(\cL_*)$ is compatible with the coaugmentation maps $\eta$.
Under these circumstances it is proven in [Br-V-(12.11)] that the induced map 
$\Gamma_c(U, \cL_*) \to \Gamma_c(U, \cD\cD(\cL_*))$
is a quasi-isomorphism.
\end{proof}

In the following theorem we use the notion of cohomological dimension $\dim_\ZZ X$ as given
in [Br-II, \S 16].
\bigskip 

\noindent (3.7) {\bf Theorem.}([Br-V-(12.19)])\quad{\it
Let $h:\cA_*\to \cB_*$ be a map of complexes of $c$-soft sheaves on $X$, and assume it induces 
a quasi-isomorphism 
$$h: \Gamma_c(U, \cA_*)\to \Gamma_c(U, \cB_*)$$
for each open $U$. 
If  $\Phi$ paracompactifying and $\dim_\ZZ X<\infty$, then for each $U$ the induced map 
$$h: \Gamma_\Phi(U, \cA_*) \to \Gamma_\Phi(U, \cB_*)$$
is a quasi-isomorphism.
}\bigskip 

The next theorem tells us that a quasi-resolution arises from the complex $\cS_*$. 
\bigskip 

\noindent (3.8) {\bf Theorem.}([Br-V-(12.14)])\quad{\it
If $X$ is locally contractible, then $\fS_*=\Gamma_c\{\cS_*\}$ is a quasi-coresolution of $\ZZ$. 
}\bigskip 

The proof is straightforward from the definitions. 
Note that our space $X$ is  locally contractible by assumption, therefore 
$\Gamma_c\{\cS_*\}$ is a quasi-coresolution of $\ZZ$ by  (\fScoreslution). 
\bigskip

\noindent (3.9) {\bf Proposition.}\quad{\it 
The following diagram commutes:
$$\begin{array}{ccccc}
&& \ti{S}_*(U) &\mapr{\Theta}  & \Gamma(U, \cS_*) \\
&\raise1ex\hbox{$\scriptstyle{\xi}$}\!\!\swarrow &\mapd{\frak{d}} &     &\mapdr{\frak{d}} \\
D(\chS^*(U))&\mapl{}&\widetilde{DDS}_*(U)   &\mapr{\Theta''}  &\ph{\,.}\Gamma(U, \cD\cD(\cS_*))\,.
 \end{array}
$$  
If we further assume $\dim_\ZZ X<\infty$, 
the maps in the above commutative diagram are all quasi-isomorphisms.
}\bigskip 

\begin{proof}  We already know the commutativity of the diagram.
We also know from (\lowerScalS), that the maps $\Theta$, $\Theta''$ 
and $\widetilde{DDS}_*\to D(\chS^*)$ are
are quasi-isomorphisms. 
By (\doubledualqisom), the map $\cS_* \to \cD\cD(\cS_*)$ induces a quasi-isomorphism 
on $\Gamma_c(U, -)$ for each open set $U$. Thus by (\cptqisPhiqis),  the map 
$$\frak{d}:\Gamma(U, \cS_* )\to \Gamma(U, \cD\cD(\cS_*))$$
is a quasi-isomorphism for each $U$. 
Since the diagram commutes it follows that the $\frak{d}$ in the middle, and consequently 
$\xi$ also, are quasi-isomorphisms. 
\end{proof}
\bigskip

{\bf \S 4. Comparison of the cap product in sheaf and singular theories}
\bigskip 

We shall compare the sheaf theoretic cap product and the singular cap product (in case $Z=X$). 
The latter is defined as follows. 
Let  
$$\cap: \ti{S}_*(X) \ts S^*(X) \to \ti{S}_*(X) \eqno{}$$
be the map which sends $\al \ts u\in  \ti{S}_m(X) \ts S^p(X)$ with 
$\al=\sum a_\sigma \sigma \in \tiS_m(X)$, $u\in S^p(X)$ to 
$$\al \cap u=\sum a_\sigma (\sigma\cap u)\,$$
(here $\sigma\cap u$ is as defined in \S 1). 
For a closed element $u\in S^p(X)$, we obtain a map of complexes
$$(-)\cap u: \tiS_*(X) \to \tiS_{*-p}(X) \,.\eqno{}$$
Here $\tiS_{*-p}(X)$ denotes the complex $\tiS_{*}(X)[p]$, obtained 
from $\tiS_{*}(X)$ by applying the shifting operation $[p]$.  Recall that by convention the shift 
$K[1]$ of a complex has differential $-d_K$. Thus the complex $\tiS_{*-p}(X)$ has
the group $\tiS_{m-p}(X)$ in homological degree 
$m$, and has differential $(-1)^p\partial$.

The induced map on homology
$$H^{lf}_m(X) \to H^{lf}_{m-p}(X) $$
depends only on the cohomology  class $[u]\in H^p(X)$, and is denoted by
$(-)\cap [u]$. 
\bigskip

Via the map $\theta$ in \S 2, we have an element $\theta(u)\in\Gamma(X, \cS^p)$, thus 
the cap product map $\cD(\cS^*)\ts \cS^* \to \cD(\cS^*)$  induces a map of complexes
$\cap \theta(u): \cD(\cS^*) \to\cD(\cS^{*-p})$. Here 
$\cD(\cS^{*-p}):= \left(\cD(\cS^{*})\right)[p]$. 
Our problem is to compare the singular cap product $\cap u$ and the sheaf-theoretic cap product
$\cap \theta(u)$.
\bigskip 

Let $S^p_{cpt}(X)\subset S^p(X)$ be the subgroup of cochains $u$ such that  $\theta(u)\in \Gamma(X, \cS^p)$ 
has compact support (equivalently, such that there exists a compact set $K$ of $U$ and an open covering $\{U_\al\}$ of 
$X-K$ such that $u|_{U_\al}=0$ for each $\al$). 
By [Br-I-(6.2)], there exists an exact sequence 
$$0\to S^p_0(X) \to S^p_ {cpt}(X) \mapr{\theta} \Gamma_c(X, \cS^p)\to 0\eqno{}\,, $$
with $S^p_0(X)$ the group as in \S 2.  Since $S_0^*(X)$ is acyclic, $\theta: S^*_ {cpt}(X)\to \Gamma_c(X, \cS^*)$ is a quasi-isomorphism.
Obviously one has $\chS^p(X)\subset S^p_{cpt}(X) $. 
\bigskip 

\noindent (4.1) {\bf Proposition.}\quad{\it The inclusion $\check{S}^*(X) \injto S^*_{cpt}(X)$ is a quasi-isomorphism. 
}
\smallskip 

\begin{proof}  Letting $\chS^p_0(X)= S^p_0(X)\cap  S^p_{cpt}(X)$, we have 
a commutative square of inclusions of complexes
$$\begin{array}{ccc}
{S_0^*}(X)    &\injto   &S_{cpt}^*(X)  \\
\mapu{} &     &\mapur{} \\
{\chS_0}^*(X)  &\injto   &\ph{\,.}\chS^*(X) \,.
 \end{array}
$$  

One has $S_{cpt}^*(X) ={S_0^*}(X)+\chS^*(X)$. 
Indeed, for $u\in S_{cpt}^p(X)$, let $K=|\theta(u)|$, and take a compact neighborhood 
$K'$ of $K$.  Let $u'\in \chS^p(X)$ be the element given by
$u'(\sigma)=u(\sigma)$ if $|\sigma|\subset \ctop{K'}$, and $u'(\sigma)=0$ otherwise. 
Then $u-u'\in {S_0^p}(X)$, since $(u-u')(\sigma)=0$ if $|\sigma|\subset \ctop{K'}$. 

Consequently, ${S_0^*}(X)/{\chS_0}^*(X)\cong S_{cpt}^*(X)/\chS^*(X)$. 
As we already know that  ${S_0^*}(X)$ is acyclic, it is enough to show that $\chS_0^*(X)$ is also acyclic.

The proof of acyclicity of $\chS_0^*(X)$ is similar to that for ${S_0^*}(X)$.  
For any open covering $\cU$ of $X$, the inclusion of the $\cU$-based singular chains 
$S_*^\cU(X)$ in $S_*(X)$ is a quasi-isomorphism (excision theorem); thus its dual  
$S^*(X)\surjto S_\cU^*(X)$ is also a quasi-isomorphism, so 
its kernel $K_\cU^*$ is an acyclic complex. 
It follows that $S_0^*(X)=\varinjlim_{\cU} K_\cU^*$ is also acyclic. 

To generalize this, for $A$ a subset of $X$,  let
$S_*^\cU(X, A)=S_*^\cU(X)/(S_*^\cU(X)\cap S_*(A))$. 
Since $S_*^\cU(X)\cap S_*(A)=S_*^{A\cap \cU}(A)$ is quasi-isomorphic to 
$S_*(A)$, the inclusion $S_*^\cU(X, A) \injto S_*(X, A)$ is a quasi-isomorphism, 
so its dual $ S^*(X, A)\surjto S^*_\cU(X, A)$
is also a quasi-isomorphism. 
Thus its kernel, which equals $K_\cU^*(X)\cap S^*(X, A)$, is acyclic. 
Passing to the limit over $\cU$, it follows that $S_0^*(X)\cap S^*(X, A)$ is also acyclic. 

Taking $A=X-K$ with $K\in\cpt(X)$ and then taking the limit over $K$, we conclude that 
$S_0^*(X)\cap \chS^*(X)$ is acyclic. 
\end{proof}
\bigskip 

Writing $\theta$ for the composition of the maps 
$\check{S}^*(X) \injto S^*_{cpt}(X)\mapr{\theta} \Gamma_c(X, \cS^*)$, we have thus:
\bigskip  

\noindent (4.2) {\bf Proposition.}\quad{\it 
The map 
$\theta: \chS^*(X) \to \Gamma_c(X, \cS^*)$
is a quasi-isomorphism; it induces an isomorphism
$H^p(\chS^*(X))\cong H^p\Gamma_c(X, \cS^*)$, which is identified with 
the cohomology with compact support $H^p_c(X; \ZZ)$.

The dual of $\theta'$, 
$$\theta': \Gamma(X, \cD(\cS^*))=D\Gamma_c(X, \cS^*) \to D(\chS^*(X))$$
is also a quasi-isomorphism.
}\bigskip

With (\lowerSbiduality), we have:
\bigskip

\noindent (\ThmlfandBM) {\bf Theorem.}\quad{\it 
One has quasi-isomorphisms 
$$ \Gamma(X, \cD(\cS^*)) \mapr{\theta'} D\chS^*(X)\mapl{\xi} \tiS_*(X)\,.$$
They give an isomorphism in the derived category $\Gamma(X, \cD(\cS^*))\isoto  \tiS_*(X)$;
they induce isomorphisms on homology, 
$$H_m(X) \mapr{\theta'} H_m(D\chS^*(X))\mapl{\xi} H^{lf}_m(X)\,. $$
}
\smallskip

The isomorphism thus obtained, $\operatorname{cano}: H_m(X)\to H^{lf}_m(X)$ is referred to as the 
{\it canonical isomorphism}.   
It can be shown to coincide with the isomorphism given in [Br-V], see (\camparisonofisomorphisms).
\bigskip

Each of the complexes $\Gamma(X, \cD(\cS^*))$, $D\chS^*(X)$, and $\tiS_*(X)$ appearing in the 
statement of (\ThmlfandBM) is covariantly functorial in $X$ for proper maps. 
It is immediate to verify that the maps $\theta'$ and $\xi$ are natural transformations between the
functors. 

We examine compatibility of these maps with cap product on both ends.
For this we introduce cap product for the middle term. 
The cup product for $S^*(X)$ gives by restriction the map
$\cup:{S}^*(X)\ts  \chS^*(X)  \to  \check{S}^*(X)$;
this in turn induces a map
$$\cap: D\check{S}^*(X) \ts S^*(X) \to  D\check{S}^*(X) $$
defined as follows. 
For $f\in D\check{S}^*(X)$ of degree $m$ and 
and $u\in \chS^p(X)$, $v\in S^{*}(X)$, let
$$\langle f\cap u, v \rangle =(-1)^{mp}\langle f, u\cup v\rangle \,.$$
Since this is parallel to the definition in \S 2, we have again the identity
$d(f\cap u)=(-1)^p df\cap u+ f\cap du$. 

For  $u\in S^p(X)$ closed, we have  a map of complexes
$$(-)\cap u: D\chS^*(X) \to D\chS^{*-p}(X) \,,\eqno{}$$
where $D\chS^{*-p}(X)$ is short for the complex $\left(D\chS^{*}(X)\right)[p]$. 
The induced map on homology depends only on the class $[u]$, and 
thus written 
$(-)\cap [u]  :H_m(D\chS^*(X)) \to H_{m-p}(D\chS^{*}(X))$. 
\smallskip 

{\bf Remark.}\quad Although it will not be used in the sequel, it is useful to note that 
one may view the map $(-)\cap u: D\chS^*(X) \to (D\chS^{*}(X))[p]$ as the dual of 
the map of complexes $u\cup (-): \chS^{*}(X)[-p]\to \chS^*(X)$. 
For this, if $K$ is a complex of abelian groups, $p$ an integer, verify that there is an isomorphism 
of complexes $D(K[-p])\to D(K)[p]$, which is given by multiplication by $(-1)^{np}$ on 
the degree $n$ part $D(K[-p])^n=\Hom(K, I^\bullet)=(D(K)[p])^n$. 
Applying the functor $D$ to $u\cup(-)$ and composing with this isomorphism gives the map
$(-)\cap u$. 
\bigskip

The map $\xi$ is also compatible with cap product:
\bigskip 

\noindent (4.4) {\bf Proposition.}\quad
{\it We have
$$\xi(\al)\cap u = \xi(\al\cap u)\,, $$
namely the following square commutes:
$$\begin{array}{ccc}
D\chS^*(X)\ts S^*(X)  &\mapr{\cap }  &D\chS^*(X) \\
\mapu{\xi\ts 1} &     &\mapur{\xi } \\
\tiS_*(X)\ts S^*(X)  &\mapr{\cap}  &\ph{\,.}\tiS_*(X)  \,.
 \end{array}
$$  
}\smallskip 

\begin{proof}
Let $\al\in\tiS_m(X)$, $u\in S^p(X)$ and $v\in \chS^{m-p}(X)$.
According to the definition of $\xi$ in \S 3, if 
$v\in S^{m-p}(X, X-K)$, 
let $\al=\al'+\al''$ with $\al'\in S_m(X)$ and $\al''\in \tiS_m(X-K)$. 
Then 
$$\langle \xi(\al\cap u), v\rangle =(-1)^{m-p} \langle v, \al'\cap u\rangle =(-1)^{m-p} (-1)^{(m-p)p} \langle u\cup v, \al'\rangle \,$$
by the definition of $\xi$ in \S 3 and by the identity $ \langle v, \al'\cap u\rangle = (-1)^{(m-p)p}\langle u\cup v, \al'\rangle $ from \S 1. 
Also, 
$$\langle \xi(\al)\cap u, v\rangle =(-1)^{mp}\langle \xi(\al), u\cup v\rangle =(-1)^{mp}(-1)^m \langle u\cup v, \al'\rangle \,. $$
Thus $\xi(\al\cap u)=\xi(\al)\cap u$. 
\end{proof}

\noindent (4.4.1) {\bf Corollary.}\quad
{\it  If $u\in S^p(X)$ is closed, we have a commutative diagram of complexes
$$\begin{array}{ccc}
 D\chS^*(X) &\mapr{\cap u}  &D\chS^{*-p}(X)  \\
 \mapu{\xi} &     &\mapur{\xi} \\
\tiS_*(X)   &\mapr{\cap u }  &\ph{\,.} \tiS_{*-p}(X) \,.
 \end{array}
$$  
The $\xi$ on the right stands for the map $\xi[p]:\tiS_{*}(X) [p]\to \left(D\chS^{*}(X)\right)[p]$. }
\bigskip

We have thus shown:
\bigskip 

\noindent (4.5) {\bf Theorem.}\quad{\it  For a closed element $u\in S^p(X)$, 
there is a commutative diagram of complexes, 
$$\begin{array}{ccc}
 \Gamma(X, \cD(\cS^*)) &\mapr{\cap \theta(u)}  &\Gamma(X, \cD(\cS^{*-p})) \\
 \mapd{\theta'} &     &\mapdr{\theta'} \\
  D\chS^*(X) &\mapr{\cap u}  &D\chS^{*-p}(X)  \\
  \mapu{\xi} &     &\mapur{\xi} \\
  \tiS_*(X)   &\mapr{\cap u }  &\ph{\,.} \tiS_{*-p}(X) \,. 
 \end{array}
$$  
Hence  the induced diagram on homology commutes:
$$\begin{array}{ccc}
 H_m(X) &\mapr{\cap [\theta(u)]}  &H_{m-p}(X) \\
 \mapd{\operatorname{cano}} &     &\mapdr{\operatorname{cano}} \\
  H^{lf}_m(X) \ph{a} &\mapr{\cap [u] }  &\ph{a\,.} H^{lf}_{m-p}(X) \,.
 \end{array}
$$  
}\bigskip 

\noindent (4.6)\quad 
We indicate how to prove that the canonical isomorphism in (\ThmlfandBM) coincides with 
the one in [Br-V].  This fact will not used in the sequel of the paper, so the reader may skip this 
paragraph.

For this we take another look at Proposition (\lowerScalS), and 
now treat the complex $DS_*(U)$ as something similar to the complex $S^*(U)$. 
To elucidate the similarities, we introduce the notation 
$T^*(U)=DS_*(U)$, so $U\mapsto T^*(U)$ is a complex of presheaves.
There is an injective quasi-isomorphism $S^*(U)\injto T^*(U)$. 
Also let $\check{T}^*(U)=\widecheck{DS}_*(U)$. Note then we have 
$$D\check{T}^*(U)=\widetilde{DDS}_*(U)\,.$$
\sskp

1. \,  By the same arguments as for $S^*(U)$ (see \S 2), we show the following. 
Let $\cT^*$ be the sheaf associated with the presheaf $T^*$. 
The presheaf is conjunctive, and the canonical map $\theta: T^*(U)\to \cT^*(U)$ is surjective; 
if $T_0^*(U)$ is defined to be the kernel of the map, it is acyclic, hence $\theta$ is a quasi-isomorphism.

One has a map of complexes $\cS^*\to \cT^*$, which is  a quasi-isomorphism, thus there is an induced 
quasi-isomorphism $\cS^*\to \cT^*$. 
The complex of sheaves $\cT^*$ is a  quasi-resolution of $\ZZ$, and it consists of flabby sheaves. 
\bskp

2.\, Also we can repeat (\CheckCptQis) for $T^*(U)$. 
We define the subcomplex $T^*_{cpt}(U)$ of $T^*(U)$ so that
there is an exact sequence
$$0\to T^p_0(U) \to T^p_ {cpt}(U) \mapr{\theta} \Gamma_c(U, \cT^p)\to 0\,.$$
Thus $\theta$ is a quasi-isomorphism.
We  have inclusion $\check{T}^*(U)\subset T^*_{cpt}(U)$, which is verified to be a quasi-isomorphism.

The map $\theta: \check{T}^*(U)\to \Gamma_c(U, \cT^*)$ induces a quasi-isomorphism
$$\theta': \Gamma(U, \cD(\cT^*))\to D(\check{T}^*(U))\,,$$
and the following diagram commutes.
$$
\begin{array}{ccc}
D(\check{T}^*(U))    &\mapr{}    &  D(\chS^*(U))    \\
\mapu{\theta'}&    & \mapur{\theta'}     \\
\Gamma(U, \cD(\cT^*))&\mapr{}   & \phantom{\,.}\Gamma(U, \cD(\cS^*))\,.
\end{array}
$$
\bskp

3.\, We give a map of complexes of sheaves $\al: \cD(\cS_*)\to \cT^*$. 
Recall from \S 2 that the map 
$\Theta: S_*(U)\to \Gamma_c(U, \cS_*)$  induces a map 
$\Theta': \Gamma(U, \cD(\cS_*))\to T^*(U)$. Composing with the canonical map 
$\theta: T^*(U)\to \cT^*(U)$, we obtain $\al(U)$.

The composition of the maps 
$$\Gamma(X, \cD(\cT^*))\mapr{\theta'} D\check{T}^*(X)\mapr{\Theta''}\Gamma(X, \cD\cD(\cS_*))$$
is coincides with the map $\Gamma(U, \cD(\al))$.  The verification is immediate and left to the reader.
\bskp

4.\, 
From these we obtain a commutative diagram (enlargening the one before (\Sectioncsoftsheaf))
$$\begin{array}{ccccc}
&& \ti{S}_*(X) &\mapr{\Theta}  & \Gamma(X, \cS_*) \\
&\raise1ex\hbox{$\scriptstyle{\xi}$}\!\!\swarrow &\mapd{\frak{d}} &     &\mapdr{\frak{d}} \\
D(\chS^*(X))&\mapl{}&\widetilde{DDS}_*(X)   &\mapr{\Theta''}  &\Gamma(X, \cD\cD(\cS_*)) \\
\mapu{\theta'}&&\mapu{\theta'}& \,\nearrow \!\!\!\lower1ex\hbox{$\scriptstyle{\cD(\al)}$}\!\! & \\
\Gamma(X, \cD(\cS^*))&\mapl{} &\phantom{\,.}\Gamma(X, \cD(\cT^*))\,.&&
 \end{array}
$$  
\bskp

5.\, Let $\cL_*$ is a complex of flabby sheaves such that $\Gamma_c\{\cL_*\}$ is a 
quasi-coresolution of $\ZZ$, and assume $\dim_Z X<\infty$. Then by [Br-V, (12.20)] there is an isomorphism
$$H_p(X)\cong H_p\Gamma(\cL_*)$$ obtained as follows. 

The complex $\cD(\cL_*)$ is a quasi-resolution of $\ZZ$ by (\Dualcoresolution); 
if $\cI^*$ is an injective resolution of $\ZZ$, there exists a map of quasi-resolutions 
$\cD(\cL_*)\to \cI^*$. 
Recall from \S 2 that there is a map 
$\frak{d}:\cL_*\to \cD\cD(\cL_*)$.  We have thus maps of complexes 
$\cL_*\to \cD\cD(\cL_*)\gets \cD(\cI^*)$. 
It is proven that the induced maps on global sections 
$$\Gamma(X, \cL_*)\to \Gamma(X, \cD\cD(\cL_*))\gets \Gamma(X, \cD(\cI^*))$$
are quasi-isomorphisms, hence the stated isomorphism on homology. 

In our situation one has a map of quasi-resolutions $\al: \cD(\cS_*)\to \cT^*$, thus in the above we may replace 
$\cI^*$ by $\cT^*$. So the quasi-isomorphisms
$\Gamma(X, \cS_*)\to \Gamma(X, \cD\cD(\cS_*))\gets \Gamma(X, \cD(\cT^*))$
give rise to the isomorphism of [Br]. 
From the commutative diagram of item 4, we see that it coincides with our canonical isomorphism.
\bskp

\newpage

{\bf \S \SectcomparisonBMsupp. Comparison of the supported cap product in sheaf and singular theories }
\bigskip

{\bf The complexes $\chS^*(U)_X$ and $\chS^*_Z(X)$.}
\quad
Unlike the complex $\Gamma_c(X, \cS^*)$, 
the complex $ \chS^*(X)$ is not covariantly functorial for inclusions of open sets; but there is  functoriality 
in the derived category. 

For an open set $U$ of $X$, let $Cpt(U)$ denote the set of compact sets $K\subset U$, and 
let 
$$\chS^p(U)_X=\varinjlim_{K\in Cpt(U)} S^p(X, X-K) \subset \chS^p(X)\,,$$
 so  one has a subcomplex $\chS^*(U)_X\subset \chS^*(X)$. 
For each $K$, the  surjection $S^*(X, X-K) \to S^*(U, U-K)$ is a quasi-isomorphism by the excision theorem, 
so the induced surjection $\chS^*(U)_X\to \chS^*(U)$ is also a quasi-isomorphism. 

For a smaller open set $V\subset U$, one has $\chS^*(V)_X\subset \chS^*(U)_X$. 
The quasi-isomorphism $\chS^*(V)_X\to \chS^*(V)$ factors as 
$\chS^*(V)_X\to \chS^*(V)_U\to \chS^*(V)$, where the first map is the quasi-isomorphism
obtained from the quasi-isomorphisms $S^*(X, X-K) \to S^*(U, U-K)$ for $K\in \cpt(V)$ by 
taking direct limit. We have a commutative diagram of complexes 
\[
  \xymatrix@R=10pt{
 \chS^*(V)     &     &    \\
\, \chS^*(V)_U\ar[u] \ar[r] &\chS^*(U)     &    \\
 \, \chS^*(V)_X \ar[r]\ar[u] &\,\chS^*(U)_X \ar[u] \ar[r]   & \chS^*(X)
  }   
\]
Hence follows the transitivity of the maps $\chS^*(U)\to \chS^*(X)$ 
in the derived category.
\bigskip

Let $Z$ be a locally contractible closed subset of $X$. 
Set $U= X-Z$. 
We define the complex $\chS_Z^*(X)$ to be  $\check{S}^*(X)/\check{S}^*(U)_X$. 
The restriction map $\chS^*(X)\to \chS^*(Z)$ induces a map 
$\chS_{Z}^*(X)\to \chS^*(Z)$.

We show that the complexes we just introduced compare well with the sheaf theoretic counterparts:
\bigskip 

\noindent (\SectcomparisonBMsupp.1) {\bf Proposition.}\quad{\it  The map $\chS_{Z}^*(X)\to \chS^*(Z)$ is a quasi-isomorphism.
We also have a commutative square of complexes
$$\begin{array}{ccc}
\chS_{Z}^*(X)    &\mapr{}  &\chS^*(Z)  \\
\mapd{\theta} &     &\mapdr{\theta} \\
 \Gamma_c(X, \cS^*)/\Gamma_c(U, \cS^*)&\mapr{}  &\Gamma_c(Z, \cS^*_Z)
 \end{array}
$$  
with $\cS_{Z}^*$ be the singular cochain complex on $Z$,  where all maps are quasi-isomorphisms. 
}\smallskip 

\begin{proof} 
We first note that the diagram of complexes
$$\begin{array}{cccccc}
 \Gamma_c(U, \cS^*)  & &\xrightarrow{\ph{aaa}j_!\ph{aaa} } & &\Gamma_c(X, \cS^*)\\
  \mapu{\theta} & &  & &\mapu{\theta} \\
\chS^*(U)   &\leftarrow  &\chS^*(U)_X  &\injto &\chS^*(X)
 \end{array}$$
commutes, where $j_!$ is extension by zero associated with inclusion $j:U\to X$.
It is enough to show that for $K\in\cpt(U)$ the square
$$\begin{array}{ccc}
\Gamma_c(U, \cS^*)    &\mapr{j_!}  &\Gamma_c(X, \cS^*)  \\
\mapu{\theta} &     &\mapur{\theta} \\
 S^*(U, U-K)  &\mapl{}  &S^*(X, X-K)
 \end{array}
$$  
commutes, namely for $u\in S^p(X, X-K)$, one has $j_!(\theta(u|_U))= \theta (u)$.
The element $\theta(u|_U)$, thus $j_!(\theta(u|_U))$ also, has support contained in $K$; 
the element $\theta (u)$ also has support in $K$.  It is thus enough to show that the restrictions to $U$ coincide, 
namely $\theta(u|_U)= \theta(u)|_U$, but this is obvious.

Therefore the left square in the diagram below commutes. 
The second row is exact since $\cS^*$ is flabby, hence $c$-soft.
Therefore there is  a unique map $\theta: \chS_{Z}^*(X)\to \Gamma_c(Z, \cS^*|_Z)$ making 
the right square commute; it is also a quasi-isomorphism.
$$\begin{array}{ccccc}
0\to \chS^*(U)_X &\mapr{}&\chS^*(X)  &\mapr{} &\chS_{Z}^*(X)\to 0  \\
\mapd{\theta} &  &\mapd{\theta}&  &\mapdr{\theta} \\
0\to\Gamma_c(U, \cS^*)&\mapr{}&\Gamma_c(X, \cS^*) &\mapr{} &\ph{\,.}\Gamma_c(Z, \cS^*|_Z)\to 0\,.
  \end{array} 
 $$
We consider the commutative diagram of complexes
$$\begin{array}{ccc}
\chS_{Z}^*(X)    &\mapr{}  &\chS^*(Z)  \\
\mapd{\theta} &     &\mapdr{\theta} \\
 \Gamma_c(Z, \cS^*|_Z)  &\mapr{}  &\Gamma_c(Z, \cS^*_Z)
 \end{array}
$$  
The vertical maps $\theta$ are quasi-isomorphisms.
The lower map is a quasi-isomorphism, since it is induced from the map $\cS^*|_Z\to \cS_Z^*$, which is a map of 
$c$-soft resolutions of $\ZZ$ on $Z$. 
The assertion hence follows.
\end{proof}


We wish to generalize the results in the previous section to the case of supported cap product. 
For this we shall show that the cup product $ {S}^*(X)\ts \chS^*(X) \to  \check{S}^*(X)$ induce a map 
 ${S}^*(X, U)\ts \chS^*(Z) \to  \check{S}^*(X)$ in the derived category. 
By the above proposition, one may replace $\chS^*(Z)$ with $\chS_{Z}^*(X)$. 
We will also replace the complex $\check{S}^*(X)$ up to quasi-isomorphism:
\bigskip

\noindent (\chSXUnatural) {\bf Definition.}\quad  
Let 
$$\check{S}^*(X)^\natural:=\check{S}^*(X)/\check{S}_0^*(X)\,$$
where $\check{S}_0^*(X)=S_0^*(X)\cap \check{S}^*(X)$ as in the proof of (\CheckCptQis).
Since the subcomplex $\check{S}_0^*(X)$ is acyclic,
the map $\check{S}^*(X)\to \check{S}^*(X)^\natural$ is a quasi-isomorphism.
Its dual $D(\check{S}^*(X)^\natural)\to D(\check{S}^*(X))$ is also a quasi-isomorphism.
\bigskip 

Consider now the restriction of the cup product $ {S}^*(X, U)\ts \chS^*(X) \to  \check{S}^*(X)$. 
\bigskip

\noindent (\suppcapprod) {\bf Proposition.}\quad{\it 
This map induces a map of complexes 
$$ {S}^*(X, U)\ts  \chS_{Z}^*(X)\to  \check{S}^*(X)^\natural\,.$$
}
\begin{proof}
Take any  elements $u\in  {S}^*(X, U)$ and $v\in \chS^*(U)_X$. 
Then $v\in S^{*}(X, X-K)$ for a compact set $K\subset U$; thus $u\cup v$ is zero on $X-K$ and 
on $U$, and $\{X-K, U\}$ covers $X$, so $u\cup v\in \check{S}_0^*(X)$. 
The assertion hence follows. 
Note this proof shows that $u\cup v$ need not be zero in $\check{S}^*(X)$.
\end{proof}

Following the pattern in \S 4, we can define the cap product
$$\cap:  D(\check{S}^*(X)^\natural)\ts S^*(X, U) \to  D(\check{S}_Z^*(X) \,)$$
by the same formula, including the sign.
If $u$ be an element of $ S^p(X, U)$ with $du=0$, then 
we obtain a map of complexes
$$(-)\cap u: D(\chS^*(X)^\natural) \to D(\chS_Z^{*-p}(X))\,. $$
One has the induced map on homology, which depends only  on the class $[u]$, 
$(-)\cap [u]  :H_m(D\chS^*(X)) \to H_{m-p}\left(D\chS_Z^{*}(X)\right)$.

We now proceed to define the supported cap product for singular homology.
\bigskip

{\bf Supported cap product for locally finite singular homology}\quad

For an element $u\in S^p(X, U)$ with $du=0$, we obtain a diagram of complexes
$$(\star)\qquad 
\begin{array}{ccccccc}
 D(\check{S}^*(X)) &\leftarrow &D(\check{S}^*(X)^\natural)   &\mapr{\cap u}  &D(\check{S}_Z^{*-p}(X)) 
  &\leftarrow
& D(\check{S}^{*-p}(Z)) \\
\mapu{\xi}&&&     &&&\mapur{\xi} \\
\ti{S}_*(X)  & &&&&&\phantom{\,.}\ti{S}_{*-p}(Z)\,.
 \end{array}
$$  
All the maps except $\cap u$ are quasi-isomorphisms. 
Inverting some of the quasi-isomorphisms, we obtain maps 
in the derived category $D(\check{S}^*(X))\to D(\check{S}^{*-p}(Z))$ and 
$\ti{S}_*(X)\to \ti{S}_{*-p}(Z)$, both of which are still written $\cap u$, 
making the following diagram commute
$$
\begin{array}{ccc}
 D(\check{S}^*(X))   &\mapr{\cap u}    &  D(\check{S}^{*-p}(Z))    \\
\mapu{\xi}&    & \mapu{\xi}     \\
\ti{S}_*(X)&\mapr{\cap u}   & \ti{S}_{*-p}(Z)\,.
\end{array}
$$
Note that the map $\cap u:\ti{S}_*(X)\to \ti{S}_{*-p}(Z)$ is obtained only through 
the map $\cap u:D(\check{S}^*(X))\to D(\check{S}^{*-p}(Z))$, and only in the 
derived category.

So we have the induced map $\cap [u]: H^{lf}_m(X) \to H^{lf}_{m-p}(Z)$, 
and a commutative diagram
$$
\begin{array}{ccc}
 H_m\left(D\check{S}^*(X)\right)  &\mapr{\cap [u]}    & H_{m-p}\left( D\check{S}^{*}(Z)\right)  \\
\mapu{\xi}&    & \mapu{\xi}     \\
H^{lf}_m(X)&\mapr{\cap [u]}   & H^{lf}_{m-p}(Z)\,.
\end{array}
$$
This is the {\it supported cap product} for locally finite singular homology. 
\smallskip

The construction has compatibility with inclusions of closed sets $Z$.  
Let $Z'$ be a closed subset of $X$ containing $Z$, and let $U'=X-Z'$. 
There is a natural injection $S^*(X, U)\to S^*(X, U')$
and a natural surjection $\chS^*_{Z'}(X) \to \chS^*_{Z}(X)$.
The cup products $ {S}^*(X, U)\ts\chS_{Z}^*(X)  \to  \check{S}^*(X)^\natural$
and ${S}^*(X, U')\ts \chS_{Z'}^*(X)\to  \check{S}^*(X)^\natural$
are compatible via these maps, namely the diagram 
$$
\begin{array}{cccccc}
{S}^*(X, U)&\ts&\chS_{Z}^*(X) &    &\mapr{}    &  \check{S}^*(X)^\natural     \\
\mapd{}&{}&\mapu{}&&    & \Vert     \\
{S}^*(X, U')&\ts& \chS_{Z'}^*(X)&  &\mapr{}   &   \check{S}^*(X)^\natural
\end{array}
$$
``commutes" in the following sense: for $u\in S^*(X, U)$
and $v'\in \chS_{Z'}^*(X)$, 
if $u'\in S^*(X, U')$ is the image of $u$ and 
$v\in \chS_{Z}^*(X)$ is the image of $v'$, 
then $u\cup v=u'\cup v'$. 

It follows that the cap product 
$\cap :  D(\check{S}^*(X)^\natural)\ts S^*(X, U) \to  D\check{S}_Z^*(X) $
and $\cap :  D(\check{S}^*(X)^\natural)\ts S^*(X, U') \to  D\check{S}_{Z'}^*(X) $
are are also compatible.

From this we deduce that, if $u\in S^p(X, U)$ is a closed element,
 which may be viewed as a closed element of 
$S^p(X, U')$, the digram $(\star)$ for $U$ and the digram $(\star)$ for $U'$
are compatible, and consequently, 
we have a commutative diagram in the derived category
$$\begin{array}{ccc}
\tiS_*(X)    &\mapr{\cap u}  & \tiS_{*-p}(Z) \\
 \Vert &     &\mapdr{} \\
\tiS_*(X)    &\mapr{\cap u}  & \phantom{\, ,}\tiS_{*-p}(Z')\,.
 \end{array}
$$ 
In particular it induces a commutative diagram on homology
$$\begin{array}{ccc}
 H^{lf}_m(X) &\mapr{\cap [u]}  & H^{lf}_{m-p}(Z) \\
 \Vert &     &\mapdr{} \\
 H^{lf}_m(X) &\mapr{\cap [u]}  & H^{lf}_{m-p}(Z')
 \end{array}
$$ 
\bigskip 

\noindent (\comparisonprep)\quad
We now relate the upper sequence of maps in the diagram $(\star)$ to the sheaf theoretic cap product. 
\smallskip 

1.  The map $\theta: \chS^*(X)\to \Gamma_c(X, \cS^*)$ induces an isomorphism
$\theta: \chS^*(X)^\natural\to \Gamma_c(X, \cS^*)$. 
\smallskip 

2.  One has  a quasi-isomorphism $\theta:\chS_{Z}^*(X)\to \Gamma_c(X, \cS^*) / \Gamma_c(U, \cS^*)$.
The map $\theta(u)\cup: \Gamma_c(X, \cS^{*-p})\to \Gamma_c(X, \cS^*)$ factors through the quotient
$\Gamma_c(X, \cS^{*-p})/\Gamma_c(U, \cS^{*-p})$. 
Thus we have a commutative diagram 
$$\begin{array}{ccccc}
&&\Gamma_c(X, \cS^*)   &\mapl{ \theta(u)\cup}  & \Gamma_c(X, \cS^{*-p})/\Gamma_c(U, \cS^{*-p})   \\
&\raise1ex\hbox{$\scriptstyle{\theta}$}\!\!\nearrow &\mapu{\theta} &     &\mapur{\theta} \\
\check{S}^*(X)&\to &\check{S}^*(X)^\natural   &\mapl{ u\cup}  &\ph{\,.} \check{S}_Z^{*-p}(X)\,.
 \end{array}
$$  
The dual of $\Gamma_c(X, \cS^{*-p})/\Gamma_c(U, \cS^{*-p})$ obviously equals 
$\Gamma_Z(X, \cD(\cS^{*-p}))$;
the dual of  $\theta:\chS_{Z}^*(X) \to \Gamma_c(X, \cS^*) / \Gamma_c(U, \cS^*)$ is a quasi-isomorphism
$\theta': \Gamma_Z(X, \cD(\cS^*)) \to D(\chS_{Z}^*(X))$.
Thus the above diagram gives us a commutative diagram
$$\begin{array}{ccccc}
&& \Gamma(X, \cD(\cS^*))    &\mapr{\cap \theta(u)}  &\Gamma_Z(X, \cD(\cS^{*-p}))      \\
&\raise1ex\hbox{$\scriptstyle{\theta'}$}\!\!\swarrow &\mapd{\theta'} &     &\mapdr{\theta'} \\
D(\chS^*(X))&\mapl{}&D(\chS^*(X)^\natural)  &\mapr{\cap u}  &\ph{\,.}D(\chS_Z^{*-p}(X))\,.
 \end{array}
$$ 

3.  From (\chSXUtochSZ), we have a commutative square
$$\begin{array}{ccc}
 \Gamma_c(X, \cS^*)/\Gamma_c(U, \cS^*)&\mapr{}  &\Gamma_c(Z, \cS^*_Z)\\
 \mapu{\theta} &     &\mapur{\theta} \\
\chS_{Z}^*(X)   &\mapr{}  &\chS^*(Z) 
\end{array}
$$  
hence also another commutative square
$$\begin{array}{ccc}
\Gamma_Z(X, \cD(\cS^{*-p}))&\mapl{} &\Gamma(Z, \cD(\cS_Z^{*-p}))\\
\mapd{\theta'} &     &\mapdr{\theta'} \\
D(\check{S}_Z^{*-p}(X))& \mapl{} &D(\check{S}^{*-p}(Z) )
\end{array}
$$  
with all maps quasi-isomorphisms. 
\smallskip 
We have established the following theorem.
\bigskip 

\noindent (\supportedSheafsingularcomparison) {\bf Theorem. }\quad{\it  There is a  commutative diagram of complexes
$$(\star\star)_X\qquad\quad
\begin{array}{ccccccc}
 &&\Gamma(X, \cD(\cS^*))  &\mapr{\cap \theta(u)}  &\Gamma_Z(X, \cD(\cS^{*-p})) &\mapl{}&\Gamma(Z, \cD(\cS_Z^{*-p})) \\
&\raise1ex\hbox{$\scriptstyle{\theta'}$}\!\!\swarrow & \mapd{\theta'} &     &\mapdr{\theta'} &&\mapdr{\theta'}\\
D(\check{S}^*(X))&\leftarrow& D(\check{S}^*(X)^\natural)   &\mapr{\cap u}  &D(\check{S}_Z^{*-p}(X))
&\mapl{}&D(\check{S}^{*-p}(Z) )  \\
\mapu{\xi}&&&     &&&\mapur{\xi} \\
 \ti{S}_*(X)  &&&  && &\ph{\,,}\ti{S}_{*-p}(Z)\,, 
 \end{array}
$$  
in which all the arrows, except the two horizontal maps $\cap u$, are quasi-isomorphisms. 
Hence we have a commutative diagram in the derived category,
$$\begin{array}{ccc}
 \Gamma(X, \cD(\cS^*))  &\mapr{\cap \theta(u)} &\Gamma(Z, \cD(\cS_Z^{*-p})) \\
 \mapd{\theta'} &              &\mapdr{\theta'}\\
 D(\check{S}^*(X))   &\mapr{\cap u} &D(\check{S}^{*-p}(Z) )\\
\mapu{\xi} &      &\mapur{\xi} \\
 \ti{S}_*(X)   &\mapr{\cap u} &\ph{\,.}\ti{S}_{*-p}(Z)\,. 
 \end{array}
$$  
It induces a commutative diagram on homology
$$\begin{array}{ccc}
 H_m(X) &\mapr{\cap [\theta(u)]}  &H_{m-p}(Z)  \\
 \mapd{\operatorname{cano}} &     &\mapdr{\operatorname{cano}} \\
 H^{lf}_m(X)  &\mapr{\cap [u]} &\ph{\,.}H^{lf}_{m-p}(Z)\,.
 \end{array}
$$  
with vertical maps the canonical isomorphisms.
}\bigskip 

Thus the sheaf theoretic supported cap product coincides with 
the supported cap product for singular homology.
\bigskip

{\bf \S \Sectsimplicial. Simplicial supported cap product; comparison to the singular cap product }
\bigskip 

Let $X$ be an abstract simplicial complex assumed  locally finite, countable, and  of finite dimension. 
The geometric realization $|X|$ will be also written $X$; it satisfies the assumption (**) in \S 2. 
\bigskip

\noindent (\simpchains)\quad
By $C_*(X)$ we denote the complex of {\it ordered} simplicial chains, and $C^*(X)$ its dual
(see [Mu, p. 76], [Sp. p.170]). 
Recall that $C_m(X)$ is the free abelian group generated by $(v_0, \cdots, v_m)$, 
where $v_0, \cdots, v_m$ are vertices of $X$ (repetition allowed) spanning a simplex of dimension $\le m$.
Let $\tiC_*(X)$ be the complex of (locally finite) infinite ordered simplicial chains. 
We let $\chC^*(X)\subset C^*(X)$ denote the subcomplex of cochains $u$ satisfying the following 
condition: there exists a finite subcomplex outside which $u$ vanishes.

The homology of the complex $C_*(X)$ (resp. $\tiC_*(X)$) is the simplicial homology (resp. locally finite simplicial homology)
of $X$.  When convenient, we write $H_m^\simp(X)$ for $H_m(\tiC_*(X))$.
\bigskip

1.\, There are natural maps of complexes, which are known to be quasi-isomorphisms,
 $C_*(X)\to S_*(X)$, $\tiC_*(X)\to \tiS_*(X)$ and 
$S^*(X)\to C^*(X)$, $\chS^*(X)\to \chC^*(X)$. 
\smallskip 

2.\, For $Y$ a subcomplex of $X$, one has relative chain and cochain complexes 
$C_*(X, Y)$, $\tiC_*(X, Y)$ and $C^*(X, Y)$ and $\chC^*(X, Y)$. 
There are natural quasi-isomorphisms between relative complexes, e.g., 
 $C_*(X, Y)\to S_*(X, Y)$ and $\tiC_*(X, Y)\to \tiS_*(X, Y)$. 
 \smallskip 
 
3.\, One has cup and cap product on the simplicial chain and cochain complexes.

\quad $\bullet$\quad Cup product $\cup: C^*(X)\ts C^*(X)\to C^*(X)$, which is compatible with the singular cup product 
via the map $S^*(X)\to C^*(X)$. 
This restricts to $C^*(X)\ts \chC^*(X)\to \chC^*(X)$, to be used later.

\quad $\bullet$\quad  Cap product $\cap: C_*(X)\ts C^*(X)\to C_*(X)$, which is compatible with singular cap product.
This extends to $\tiC_*(X)\ts C^*(X) \to \tiC_*(X)$.
\smallskip 

\noindent Thus a  cocycle $u\in C^p(X)$ induces maps of complexes
$\cap u: \tiC_*(X) \to \tiC_{*-p}(X)\,, $ and 
$ u\cup  :\chC^{*-p}(X)\to \chC^*(X)$. 
\bigskip 

\noindent (\suppcapcupprod) {\bf Supported cap and cup product for simplicial homology}
\smallskip 

Let $Z$ be a subcomplex of $X$, and $U=X-Z$ be its complement.  
Let $N=N(Z)$ be the closed star of $Z$, namely the subcomplex consisting of simplices
meeting $Z$; let $N^c$ be the subcomplex consisting of simplices disjoint from 
$Z$.
We assume that
 the inclusion $Z\subset N$ is a proper deformation retract.
 (This can be so arranged by taking barycentric subdivision of $X$ twice -- see [Mu, Lemma 70.1] for
 an argument of this kind).
\bigskip

1.\, One has cap product 
$\tiC_*(X)\ts C^*(X, N^c)\to \tiC_*(N)$. 
\begin{proof} For $u\in C^p(X, N^c)$ and an $m$-simplex $\sigma$, one has by definition
$\sigma\cap u =u(\sigma')\sigma''$. 
If the front face $\sigma'$ is a simplex of $N^c$, then $u(\sigma')=0$.
Otherwise $\sigma'$ meets $Z$, so $\sigma$ is a
simplex of $N$. In particular $\sigma''$ is a simplex of $N$.
Thus $u(\sigma')\sigma''$ is a chain in $N$.
\end{proof}

Therefore, a cocycle $u\in C^p(X, N^c)$ induces a map of complexes
$$\cap u: \tiC_*(X) \to \tiC_{*-p}(N)\,. $$
The induced map on homology
$$\cap u: H_m(\tiC_*(X)) \to H_{m-p}(\tiC_{*-p}(N)) \,,$$
which depends only on the class $[u]$, is the {\it supported cap product}
on simplicial homology. 
The inclusion $Z\to N$, being a proper deformation retract, induces an isomorphism 
$H_{m-p}(\tiC_{*}(Z))\to H_{m-p}(\tiC_{*}(N))$; the composition
$H_m(  \tiC_*(X)) \mapr{\cap [u]} H_{m-p}(\tiC_*(N))\mapl{\cong} H_{m-p}(\tiC_*(Z))$
will be also called the supported cap product, and written $\cap [u]$.

In order to compare this with the singular cap product, as preliminaries we 
develop some parallels of \S 5 for the simplicial theory.  
\smallskip

2.\, One has cup product 
$C^*(X, N^c)\ts \chC^*(N)\to \chC^*(X)$. 

\begin{proof}
There is an exact sequence of complexes $0\to \chC^*(X, N)\to \chC^*(X) \to  \chC^*(N)\to 0$. 
If $u\in C^p(X, N^c)$, $v\in  \chC^{m-p}(X, N)$ and $\sigma$ is an $m$-simplex, then 
$(u\cup v)(\sigma)=\pm u(\sigma')v(\sigma'')$, which is zero since $\sigma'\subset N^c$ or 
$\sigma''\subset N$.
\end{proof}

Therefore, if $u\in C^p(X, N^c)$ is a cocycle, one has a map of complexes
$$ u\cup: \chC^{*-p}(N)\to \chC^*(X)\,.$$ 
\smallskip 

3.\, The cup product in item 2 induces a map 
$$\cap: D(\chC^*(X))\ts C^*(X, N^c)\to D(\chC^*(N))\,,$$
as in \S 4 (and \S 5).  Hence, for a cocycle $u\in C^p(X, N^c)$
we have the map 
$\cap u:D(\chC^*(X))\to D(\chC^{*-p}(N))$. 
\smallskip 

4.\, We define a map $\xi: \tiC_*(X) \to  \Hom(\chC^*(X), \ZZ)\subset
D( \chC^*(X)) $, analogous to the map 
$\xi$ in \S 3 for $\tiS_*(X)$. 

For $\al\in \tiC_m(X)$ and $u\in \chC^m(X)$, let $Y$ be a finite subcomplex outside which $u$ 
vanishes. 
Take a finite subcomplex $Z$ containing the star of $Y$. 
We can take a decomposition $\al=\al'+\al''$, where $\al'\in C_m(Z)$, and $\al''=\sum b_\sigma \sigma\in \tiC_m(X)$ 
satisfy the condition that if $b_\sigma\neq 0$, then $\sigma$ is disjoint from $Y$. 
For example, for $\al= \sum a_\sigma \sigma$, one may take
$$\al'=\sum_{\sigma \subset Z} a_\sigma \sigma, \qquad \al''=\sum_{\sigma \not\subset Z} a_\sigma \sigma\,.$$
Let $\xi(\al)\in \Hom( \chC^m(X), \ZZ)$ be the element given by
$$\langle \xi(\al), u\rangle =(-1)^m \langle u, \al'\rangle \,.$$ 
One verifies easily that this is independent of the choice of $Z$
and of the decomposition of $\al$;
also $\xi$ is a map of complexes.
\bigskip

As we did for the map $\xi$ for $\tiS_*(X)$ (see (\xiforS) and its corollary), one shows:
\bigskip

\noindent (\xiCDC) {\bf Proposition.}\quad{\it For a cocycle $u\in C^p(X, N^c)$, the following square commutes:
$$\begin{array}{ccc}
  D( \chC^*(X))  &\mapr{\cap u}  &D(\chC^{*-p}(N)  )  \\
\mapu{\xi} &     &\mapur{\xi} \\
 \tiC_*(X)  &\mapr{\cap u}  &\ph{\,.}\tiC_{*-p}(N) \,.
 \end{array}
$$
}\bigskip

We also have, obviously from the definitions of the maps $\xi$, 
\bigskip

\noindent (\xiCS) {\bf Proposition.}\quad{\it 
The following square commutes, in which the horizontal maps are  the natural ones.
$$\begin{array}{ccc}
  D( \chC^*(X))  &\mapr{}  &D( \chS^*(X))  \\
\mapu{\xi} &     &\mapur{\xi} \\
 \tiC_*(X)  &\mapr{}  &\tiS_{*}(X) \,.
 \end{array}
$$
}\bigskip

With these preparations, we have come to the main part of this section.
 For the rest of this section, we assume given a cocycle $u\in S^p(X, X-Z)$. 
Note it induces a cocycle $\rho(u)\in C^p(X, N^c)$; we often just write 
$u$ for $\rho(u)$. 
\bigskip

The cocycle $u$ gives a map $u\cup:\check{S}_N^{*-p}(X)\to \check{S}^*(X)^\natural$, and 
$\rho(u)$ gives $\rho(u)\cup: \chC^{*-p}(N)\to \chC^*(X) $. 
So
we have a diagram of complexes
$$\begin{array}{ccccccc}
\check{S}^*(X) & \to & \check{S}^*(X)^\natural  &\mapl{u\cup}  &\check{S}_N^{*-p}(X)&\to
& \check{S}^{*-p}(N)\\
\mapd{\rho}&&&     &&&\mapdr{\rho} \\
\chC^*(X)  & &&\xleftarrow{\ph{aaaa}\rho(u)\cup\ph{aaaa} }&&&\chC^{*-p}(N)
 \end{array}
$$  
in which $\rho$ are the natural maps, and all the maps except the two maps $u\cup$ are quasi-isomorphisms. 
\bigskip

\noindent (\CScommutes) {\bf Proposition.}\quad{\it 
The above diagram induces a commutative diagram on homology.
}\bigskip

Before the proof we need to note how the map $\rho$ is affected by subdivision.
Denote by $K$ the given triangulation of the space $X$, and $K'$ its subdivision. 
The corresponding simplicial chain groups will be written $C_*(K)$ and $C_*(K')$. 
There is the subdivision map $\lambda: C_*(K)\to C_*(K')$.
If $g : K'\to K$ is a simplicial approximation of the identity map ([Mu], \S 14), one has the induced map 
$g_*: C_*(K')\to C_*(K)$, and $g_*\scirc\lambda=1$, and $\lambda\scirc g_*\simeq 1$ (homotopic to the identity).

Dually, there is the subdivision map $\lambda: C^*(K')\to C^*(K)$, the map $g^*:  C^*(K)\to C^*(K')$, and one has
$\lambda\scirc g^*=1$ and $g^*\scirc \lambda\simeq 1$.
Further, the same holds for their restrictions to subcomplexes $\lambda: \chC^*(K')\to \chC^*(K)$ and $g^*:  \chC^*(K)\to \chC^*(K')$
(see [Mu], p.269, Exercise 5). Since the map $\rho: \chS^*(X)\to \chC^*(K)$ obviously commutes with $g^*$, it follows that the 
diagram 
\[
  \xymatrix@R=10pt{
    \chC^*(K') \ar[rr]^{\lambda}  &&\chC^*(K) \\
    & \chS^*(X)\ar[ul]^{\rho} \ar[ur]_{\rho}
}
\]
induces a commutative diagram on cohomology. 
As a consequence, 
since the assertion of (\CScommutes) concerns cohomology, one may
take a subdivision of the triangulation of $X$ and replace the map $\rho$ accordingly. 

\renewcommand{\proofname}{\it{Proof of (\CScommutes)}}
\begin{proof} 
We now show the 
Let now $v\in \chS^{m-p}_N(X)$ be a cocycle, $dv=0$. 
Then one has $[u\cup v]\in H^m(\chS^*(X)^\natural)$.
Since $\chS^*(X)\to \chS^*(X)^\natural$ is a quasi-isomorphism,  there exists 
$x\in \chS^m(X)$, $dx=0$, such that $[x]\in H^m(\chS^*(X))$ maps to $[u\cup v]$.
Then there exist elements $\al\in \chS_0^m(X)$ and $\beta\in \chS^{m-1}(X)$ such that
$$x-u\cup v=\al+d\beta\,.$$

We may take a subdivision of the given triangulation, and the subdivision may be taken so that $\rho(\al)=0$.
This requires a well-known argument involving the Lebesgue number
(cf. [Mu, p. 178, Theorem 31.3])
which can be applied to 
those finite number of simplices on which $\al$ takes non-zero values.

Then one has 
$$\rho(x)-\rho(u\cup v)=d\rho(\beta)\,,$$
hence $[\rho(x)]=[\rho(u\cup v)]$, proving the assertion. 
\end{proof}
\renewcommand{\proofname}{\it{Proof.}}

For a complex $K$ and its dual $D(K)$, there is a natural short exact sequence
$$0\to \Ext^1(H^{-p+1}(K), \ZZ) \to H^p(D(K))\to \Hom(H^{-p}(K), \ZZ)\to 0$$
(see [Br, V-(2.3)]).  Thus:
\bigskip 

\noindent (\CScommutes.1) {\bf Corollary.}\quad{\it The dual of the above diagram
$$\begin{array}{ccccccc}
D(\check{S}^*(X)) & \gets & D(\check{S}^*(X)^\natural)  &\mapr{\cap u}  &D(\check{S}_N^{*-p}(X)\, )&\gets
& D(\check{S}^{*-p}(N)) \\
\mapu{\rho'}&&&     &&&\mapur{\rho'} \\
D(\chC^*(X))  & &&\xrightarrow{\ph{aaaa}\cap u\ph{aaaa} }&&&D(\chC^{*-p}(N))
 \end{array}
$$  
commutes on homology. 
}\bigskip 

%
%

We give the main result of this section:
\bigskip

\noindent (\CScapproduct) {\bf Theorem.}\quad{\it For a cocycle $u\in S^p(X, X-Z)$, 
the following diagram in the derived category
$$\begin{array}{ccc}
\tiS_*(X)    &\mapr{\cap u}  &\tiS_{*-p}(N)  \\
\mapu{\rho} &     &\mapur{\rho} \\
  \tiC_*(X)  &\mapr{\cap u }  &\tiC_{*-p}(N)
 \end{array}
$$  
induces on homology a commutative square
$$\begin{array}{ccc}
H^{lf}_m(X)    &\mapr{\cap [u]}  &H^{lf}_{m-p}(N)  \\
\mapu{\rho} &     &\mapur{\rho} \\
H^\simp_m(X)  &\mapr{\cap [u] }  &\ph{\,}H^\simp_{m-p}(N)\,.
 \end{array}
$$  
Thus the singular supported cap product coincides with the simplicial supported cap product 
via the isomorphism $H^{lf}_m(X)=H^\simp_m( X)$. 
}\smallskip 

\begin{proof} Consider the following diagram in the derived category.
$$
\xymatrix@C=15pt@R=8pt{
& D\chS^*(X) \ar[rr]^{\cap u}   &&   D(\chS^{*-p}(N))      \\
\tiS_*(X) \ar[rr]^(0.6){\cap u}  \ar[ru] &&\tiS_{*-p}(N) \ar[ru]   \\
& DC^*(X) \ar[rr]|(.47)\hole^(0.4){\cap u}  \ar[uu]|(.5)\hole  &&   D(\chC^{*-p}(N))  \ar[uu]    \\
\tiC_*(X) \ar[rr]^(0.6){\cap u}  \ar[ru] \ar[uu] &&\,\,\tiC_{*-p}(N) \ar[ru]    \ar[uu]\,\,.
 }
$$
The squares on both sides commute on chain level by (\xiCS).
The top square commutes in the derived category by our definition of $\cap u$ in \S 5.
The bottom square commutes on chain level by (\xiCDC). 
The back square induces a commutative diagram on homology, by (\CScommutes). 
It now follows that the front square induces a commutative diagram on homology.
\end{proof}

\noindent (\CScapproduct.1) {\bf Corollary.}\quad{\it The following square commutes:
$$\begin{array}{ccc} 
H^{lf}_m(X)    &\mapr{\cap [u]}  &H^{lf}_{m-p}(Z)  \\
\mapu{\rho} &     &\mapur{\rho} \\
H^\simp_m(X)  &\mapr{\cap [u] }  &\ph{\,}H^\simp_{m-p}(Z)\,.
 \end{array}
$$  
}\smallskip 

\begin{proof} The given cocycle $u\in S^p(X, X-Z)$ gives a cocycle
 in $S^p(X, X-N)$, and by 
 the compatibility of cap product with inclusions of $Z$ (see \S \SectcomparisonBMsupp) the 
diagram
$$\begin{array}{ccc}
H^{lf}_m(X)    &\mapr{\cap [u]}  & H^{lf}_{m-p}(Z) \\
 \Vert &     &\mapdr{} \\
H^{lf}_m(X) &\mapr{\cap [u]}  & H^{lf}_{m-p}(N)
 \end{array}
$$  
commutes.  
The right arrow is an isomorphism, since $Z\subset N$ is a proper deformation retract. 

There is an analogous commutative square 
$$\begin{array}{ccc}
H^\simp_m(X)   &\mapr{\cap [u]}  & H^\simp_{m-p}(Z) \\
 \Vert &     &\mapdr{} \\
H^\simp_m(X)    &\mapr{\cap [u]}  & H^\simp_{m-p}(N)
 \end{array}
$$  
by the definition of $\cap [u]$ we gave in this section.
Our assertion follows from these facts and the theorem.
\end{proof}

In particular,  the simplicial supported cap product is independent of the choice of 
a triangulation of the locally compact Hausdorff space $X$ satisfying (**) in \S 3, as long as 
$X$ has a triangulation satisfying the required conditions.
\bigskip

{\bf \S \SectBMpairs. Borel-Moore homology of pairs of spaces }
\bigskip 


Let $X$ be a locally compact Hausdorff space satisfying (**) in \S 3. 
Let $Y$ be a locally contractible closed set of $X$, and $i: Y\to X$ be the inclusion;
we will consider the homology of such a pair $(X, Y)$. 

As discussed in the proof of (\chSXUtochSZ), there is a canonical surjective map of $\cS^*|_Y=i^*\cS^*\to 
\cS_Y^*$; we examine this more closely.
\bigskip

\noindent (\StoSY) {\bf Lemma.}\quad{\it The surjection $\cS^*\to i_*\cS_Y^*$ has a canonical section. 
(The section is not a map of complexes.)
}\smallskip 

\begin{proof}  For an open set $U$ of $X$, the restriction map $S^p(U)\to S^p(U\cap Y)$ is surjective, and 
there is a section $s$ of this map given by extension by zero: For $u\in S^p(U\cap Y)$, its extension $s(u)\in S^p(U)$ is given by 
$s(u)(\sigma)=u(\sigma)$ if $\sigma\in S_p(U\cap Y)$ and zero otherwise. 

Assume in general that $P$, $Q$ are presheaves on $X$ and $Y$, respectively; assume 
that $h: P \to i_*Q$ is a map of presheaves and $s: i_*Q\to P$
is a section of it, namely $hs=id$ holds.
Letting $\ti{P}$ (resp. $\ti{Q}$) denote the associated sheaf of $P$ (resp. $Q$), one has induced maps 
of sheaves $\ti{P}\mapr{h} i_*\ti{Q}\mapr{s}\ti{P}$ with $hs=id$, thus also has 
maps 
$$i^*\ti{P}\mapr{h} \ti{Q}\mapr{s}i^*\ti{P}$$
 with $hs=id$. 
 
We apply this to the  presheaves on $X$ and $Y$, $P=S^p$ and $Q=S^p_Y$, and obtain the 
assertion. 
\end{proof}

By the lemma, any additive covariant functor applied to $\cS^* \to 
i_*\cS_Y^*$ gives a surjective map with a section.
In particular the induced map
$ \Gamma_c(X, \cS^*)\to \Gamma_c(Y, \cS^*_Y)$
is also surjective with a section, hence its dual
$$\Gamma(Y, \cD(\cS^*_Y)) \to \Gamma(X, \cD(\cS^*))$$
is an injection.  
By definition the Borel-Moore homology of the pair $(X, Y)$ is the homology of the quotient complex:
$$H_m(X, Y; \ZZ)=H_m(X, Y):=H_m\left(\Gamma(X, \cD(\cS^*))/\Gamma(Y, \cD(\cS^*_Y))\,\right)\,.$$
The reader can show that this coincides with the definition using the injective resolution in [Br, V, \S 5]. 
One has, of course, the long exact sequence of the groups
$$\to H_m(Y)\to H_m(X)\to H_m(X, Y)\to H_{m-1}(Y)\to \cdots\,.$$

Let $Z$ be another locally contractible closed set of $X$ such that $Z\cap Y$ is also locally contractible. 
Generalizing the construction in \S 2, we shall produce the supported cap product map
$$H_m(X, Y)\ts H^p_Z(X) \to H_{m-p}(Z, Z\cap Y)\,.$$

As a preliminary, let $\cL^*$ be a differential sheaf, and assume given a map of sheaves 
$\cup: \cL^*\ts \cL^*\to \cL^*$ (cup product), which is associative and satisfies
the identity $d(x\cup y)=(-1)^p dx\cup y+ x\cup dy$ if $y$ is a section of $\cL^p$. 
Then just as we described for $\cS^*$ in \S 2, one has cap product  $\cD(\cL^*)\ts \cL^*\to \cD(\cL^*)$. 

Let $\cM^*$ be another differential sheaf with cup product, and assume given a map of differential sheaves 
$h: \cL^*\to \cM^*$ that is compatible with product. 
Then we have an induced map $h_*: \cD(\cM^*)\to \cD(\cL^*)$ by functoriality. 
For sections $f$ of $\cD(\cM^*)$ and $s$ of $\cL^*$, one verifies the identity (``projection formula")
$$h_*(f)\cap s= h_*(f\cap h(s))\,.$$

Applying this to the map $\cS^*\to i_*\cS_Y^*$, by (\StoSY) we get an injective map 
$$h_*: \cD(i_*\cS_Y^*) =i_*\cD(\cS_Y^*) \to \cD(\cS^*)$$
and a diagram 
$$\begin{array}{ccccc}
 \cD(\cS^*)&\ts & \cS^* &\mapr{\cap } & \cD(\cS^*)  \\
\mapu{h_*} &  &\mapd{h}&  &\mapur{h_*} \\
 i_*\cD(\cS_Y^*)&\ts &i_* \cS_Y^*  &\mapr{\cap } &\ph{\,,} i_*\cD(\cS_Y^*) \,,
  \end{array}
 $$
which ``commutes" in the sense that the projection formula holds.  
From this we obtain a ``commutative" diagram 
$$\begin{array}{ccccc}
 \Gamma(X, \cD(\cS^*))&\ts & \Gamma_Z(X, \cS^* ) &\mapr{\cap } &\Gamma_Z(X,  \cD(\cS^*))  \\
\mapu{h_*} &  &\mapd{h}&  &\mapur{h_*} \\
 \Gamma(Y, \cD(\cS_Y^*))&\ts &\Gamma_{Z\cap Y}(Y, \cS_Y^*)  &\mapr{\cap } &\ph{\,,} \Gamma_{Z\cap Y}(Y, \cD(\cS_Y^*)) \,,
  \end{array}
 $$
 in which the two vertical maps $h_*$ are injections by (\StoSY).
So the the cap product map 
$\cap: \Gamma(X, \cD(\cS^*))\ts \Gamma_Z(X, \cS^*)\to \Gamma_Z(X, \cD(\cS^*))$
sends $\Gamma(Y, \cD(\cS_Y^*))\ts \Gamma_Z(X, \cS^*)$ into $\Gamma_{Z\cap Y}(Y, \cD(\cS_Y^*))$,
inducing a map 
$$\cap: \frac{\Gamma(X, \cD(\cS^*))}{\Gamma(Y, \cD(\cS_Y^*))}\ts \Gamma_Z(X, \cS^*)\to 
\frac{\Gamma_Z(X, \cD(\cS^*))}{\Gamma_{Z\cap Y}(Y, \cD(\cS_Y^*))}\,.$$

For the target, we have maps
$$\frac{\Gamma_Z(X, \cD(\cS^*))}{\Gamma_{Z\cap Y}(Y, \cD(\cS_Y^*))}
=\frac{\Gamma(Z, \cD(\cS^*|_Z))}{\Gamma(Z\cap Y, \cD(\cS_Y^*|_{Z\cap Y}))}
\gets \frac{\Gamma(Z, \cD(\cS^*_Z))}{\Gamma(Z\cap Y, \cD(\cS^*_{Z\cap Y}))}
\,.
$$
The first equality follows from (\PropgammaZDL).
The second map is a quasi-isomorphism, obtained as follows. 
As noted in the proof of (\chSXUtochSZ), the map $\cS^*|_Z\to \cS_Z^*$ 
induces a quasi-isomorphism
$\Gamma_c(Z, \cS^*|_Z)\to \Gamma_c(Z, \cS_Z^*)$, which induces a quasi-isomorphism
$\Gamma(Z, \cD(\cS_Z^*))\to \Gamma(Z, \cD(\cS^*|_Z))$; similarly one has a quasi-isomorphism
$\Gamma(Z\cap Y, \cD(\cS_{Z\cap Y}^*))\to \Gamma(Z, \cD(\cS_Y^*|_{Y\cap Z}))$.  

Therefore,  upon taking cohomology we have the map $H_m(X, Y)\ts H^p_Z(X) \to H_{m-p}(Z, Z\cap Y)$ as desired. 
\bigskip 

\noindent (\chSrelative)\quad
Before discussing supported cap product for singular homology, we 
introduce the relative versions of the complexes $\chS^*(U)_X$ and $\chS_Z^*(X)$.
\bigskip

1.\,  The map $\chS^*(X)\to \chS^*(Y)$ is surjective. Similarly, for each open set $U$ of $X$, the restriction map
$\chS^*(U)_X\to \chS^*(U\cap Y)_Y$ is surjective.
\smallskip

\begin{proof}
Let $u\in S^p(Y)$ be any element; for a set $L\in Cpt(Y)$, $u\in S^p(Y, Y-L)$. 
If $u'\in S^p(X)$ is its extension by zero to $X$, then one shows $u'\in S^p(X, X-L)$.
Thus the first assertion holds. The proof of the second statement is similar. 
\end{proof}

Define $\chS^*(X, Y)$ to be the kernel of the map $\chS^*(X)\to \chS^*(Y)$; similarly define
 $\chS^*(U, U\cap Y)_X$ to be the kernel of the map $\chS^*(U)_X\to \chS^*(U\cap Y)_Y$. 
\smallskip

2. \,  It follows from item 1 that the map $\chS_Z^*(X)\to \chS_{Z\cap Y}^*(Y)$ is also surjective.
Let $\chS_Z^*(X, Y)$ be the kernel of this map.  
We then deduce, using the 9-lemma,  a short exact sequence 
$$0\to \chS^*(U, U \cap Y)_X\to \chS^*(X, Y)\to \chS_{Z}^*(X, Y)\to 0\,.$$
\smallskip 

3.\, There is a natural surjective quasi-isomorphism
$\chS_{Z}^*(X, Y)\to \chS^*(Z, Z\cap Y)\,.$
\smallskip 

\begin{proof}  We have an induced map $\chS_{Z}^*(X, Y)\to \chS^*(Z, Z\cap Y)$, making a commutative diagram 
of complexes with exact rows
$$
\begin{array}{ccccccc}
0 \to &\chS_Z^*(X, Y)   &\mapr{}    &\chS_Z^*(X)     &\mapr{} &\chS_Z^*(Z\cap Y) &\to 0\\
& \mapd{}&    & \mapd{}   & &\mapd{}  &\\
0 \to &\chS^*(Z, Z\cap Y)   &\mapr{}    &\chS_Z^*(Z)     &\mapr{} &\chS^*(Z\cap Y) &\ph{\,.}\to 0\,,
\end{array}
$$
where the second and third vertical maps are quasi-isomorphism;
so the first vertical map is also a quasi-isomorphism.
\end{proof}

4.\, The restriction map $\chS^*(X)\to \chS^*(Y)$ obviously induces a 
surjective map $\chS^*(X)^\natural\to \chS^*(Y)^\natural$. 
Let $\chS^*(X, Y)^\natural$ be its kernel. 
Then the natural map $\chS^*(X, Y)\to \chS^*(X, Y)^\natural$ is a quasi-isomorphism. 
\smallskip 

5.\, Since the cup product maps 
${S}^*(X, U)\ts  \chS_{Z}^*(X)\to  \check{S}^*(X)^\natural$ and
${S}^*(Y, Y\cap U)\ts  \chS_{Z\cap Y}^*(X)\to  \check{S}^*(Y)^\natural$
are compatible via the respective restriction maps, there is an induced map 
$$\cup: {S}^*(X, U)\ts  \chS_{Z}^*(X, Y)\to  \check{S}^*(X, Y)^\natural\,.$$
This induces a cap product
$$\cap: D(\check{S}^*(X, Y)^\natural\,)\ts {S}^*(X, U)\to D(\chS_{Z}^{*}(X, Y))\,.$$
In particular, for $u\in S^p(X, U)$ closed, one has a map of complexes 
$$\cap u: D(\check{S}^*(X, Y)^\natural\,)\to  D(\chS_{Z}^{*-p}(X, Y))\,.$$
\smallskip 

6.\, One has a map of complexes $\xi: \tiS_*(X, Y)\to D(\check{S}^*(X, Y))$ which fits in a commutative diagram 
with exact rows
$$
\begin{array}{ccccccc}
 0\to &D(\chS^*(Y))  &\to    & D(\chS^*(X))   &\to  &D(\chS^*(X, Y))& \to 0   \\
 &\mapu{\xi} &&\mapu{\xi}      & &\mapur{\xi}  &  \\
0\to &\tiS_*(Y)&\to   &\tiS_*(X)  &\to     &\tiS_*(X,Y)& \to\ph{\,.} 0\,;
\end{array}
$$
thus the induced map $\xi$ is also a quasi-isomorphism.
\bigskip

{\bf Singular homology of pairs and supported cap product}\smallskip

In \S 5 we have defined the supported cap product on the locally finite singular homology, 
which we now generalize to the case of locally finite singular homology of the pair $(X, Y)$.
In that section we considered the diagram $(\star)$ for $X$ with respect to $Z$ -- which we 
name $(\star)_X$; we have also diagram  $(\star)_Y$ for $Y$ with respect to $Z\cap Y$. 
Similarly we have a diagram
$$(\star)_{X, Y}\qquad 
\begin{array}{ccccccc}
 D(\check{S}^*(X, Y)) &\leftarrow &D(\check{S}^*(X, Y)^\natural)   &\mapr{\cap u}  &D(\check{S}_Z^{*-p}(X, Y)) 
  &\leftarrow
& D(\check{S}^{*-p}(Z, Z\cap Y)) \\
\mapu{\xi}&&&     &&&\mapur{\xi} \\
\ti{S}_*(X, Y)  & &&&&&\phantom{\,.}\ti{S}_{*-p}(Z, Z\cap Y)\,.
 \end{array}
$$  
We have already explained the maps involved, which are quasi-isomorphisms except $\cap u$. 
There is a sequence of diagrams 
$$0\to (\star)_Y\to (\star)_X \to (\star)_{X, Y} \to 0$$
which is termwise exact. 

From $(\star)_{X, Y}$ we get a commutative square in the derived category
$$
\begin{array}{ccc}
 D(\check{S}^*(X, Y))   &\mapr{\cap u}    &  D(\check{S}^{*-p}(Z, Z\cap Y))    \\
\mapu{\xi}&    & \mapur{\xi}     \\
\ti{S}_*(X, Y)&\mapr{\cap u}   & \ti{S}_{*-p}(Z, Z\cap Y)\,.
\end{array}
$$
On cohomology we have the map $\cap [u]: H^{lf}_m(X, Y) \to H^{lf}_{m-p}(Z, Z\cap Y)$, as wanted.
This and the cap products $H^{lf}_m(X) \to H^{lf}_{m-p}(Z)$, $H^{lf}_m(Y) \to H^{lf}_{m-p}(Z\cap Y)$
fit in the long exact sequences for $(X, Y)$ and $(Z, Z\cap Y)$.
\bigskip

{\bf Comparison of relative Borel-Moore homologies and supported cap products}\smallskip 

Given the preparations thus far, 
the statements and the proofs of the following two results for $(X, Y)$
are parallel to the corresponding theorems for $X$. 
The first one generalizes (\ThmlfandBM).
\bigskip

\noindent (\ThmlfandBMpair) {\bf Theorem.}\quad{\it 
One has quasi-isomorphisms 
$$ {\Gamma(X, \cD(\cS^*))}/{\Gamma(Y, \cD(\cS_Y^*))}\,
\mapr{\theta'} \,D\chS^*(X, Y)\mapl{\xi} \tiS_*(X, Y)\,.$$
They induce isomorphisms on homology, 
$$\operatorname{cano}: H_m(X, Y) \to H_m(D\chS^*(X, Y))\gets H^{lf}_m(X, Y)\,.$$
 }\smallskip 

\begin{proof} 
By (\dualThetaQis) we have quasi-isomorphisms 
$\theta': \Gamma(X, \cD(\cS^*))\to D(\chS^*(X))$ and
$\theta': \Gamma(Y, \cD(\cS_Y^*))\to D(\chS^*(Y))$, 
whence the quasi-isomorphim $\theta'$ on the left between the quotients. 
The quasi-isomorphism $\xi$ on the right was given in (\chSrelative).
\end{proof}
\bigskip

The second result generalizes (\supportedSheafsingularcomparison), and the argument proceeds as follows.
\smallskip 

1.\, 
In (\comparisonprep) we have explained the quasi-isomorphism
$\theta': \Gamma_{Z}(X, \cD(\cS^*))\to D(\chS_Z^*(X))$;  a quasi-isomorphism
$\theta': \Gamma_{Z\cap Y}(Y, \cD(\cS_Y^*))\to D(\chS_{Z\cap Y}^*(Y))$ is similarly obtained. 
The second row of the following diagram is exact, as the dual of the exact sequence defining $\chS_Z^*(X, Y)$.
It follows that there is a map 
$\theta': {\Gamma_Z(X, \cD(\cS^*))}/{\Gamma_{Z\cap Y}(Y, \cD(\cS_Y^*))}\to D(\chS_Z^*(X, Y))$
which makes the whole diagram commute
and which is a quasi-isomorphism.
$$
\begin{array}{ccccccc}
 0\to &\Gamma_{Z\cap Y}(Y, \cD(\cS_Y^*))   &\mapr{}    &\Gamma_{Z}(X, \cD(\cS^*))     &\mapr{} &
 {\Gamma_{Z}(X, \cD(\cS_Y^*))}/{\Gamma_{Z\cap Y}(Y, \cD(\cS_Y^*))}      &\to 0\\
& \mapd{}&    & \mapd{}   & &\mapdr{\theta'}  &\\
0\to &D(\chS_{Z\cap Y}^*(Y))  &\mapr{}    &D(\chS_Z^*(X))   &\mapr{} &D(\chS_Z^*(X, Y))    &\to 0\,.
\end{array}
$$
\smallskip

2.\,
One has from (\comparisonprep) a commutative square
$$\begin{array}{ccc}
 \Gamma(X, \cD(\cS^*))    &\mapr{\cap \theta(u)}  &\Gamma_{Z}(Y, \cD(\cS^{*-p}))      \\
\mapd{\theta'} &     &\mapdr{\theta'} \\
D(\chS^*(X)^\natural)  &\mapr{\cap u}  &\ph{\,.}D(\chS_{Z}^{*-p}(X))\,,
 \end{array}
$$ 
another similar commutative square
$$\begin{array}{ccc}
 \Gamma(Y, \cD(\cS_Y^*))    &\mapr{\cap \theta(u)}  &\Gamma_{Z\cap Y}(Y, \cD(\cS_Y^{*-p}))      \\
\mapd{\theta'} &     &\mapdr{\theta'} \\
D(\chS^*(Y)^\natural)  &\mapr{\cap u}  &\ph{\,.}D(\chS_{Z\cap Y}^{*-p}(Y))\,,
 \end{array}
$$ 
and a termwise injection from the second to the first by (\StoSY).  Passing to the quotient, and taking item 1 into account, 
we get a commutative diagram 
$$\begin{array}{ccc}
 \Gamma(X, \cD(\cS^*))/\Gamma(Y, \cD(\cS_Y^*))   &\mapr{\cap \theta(u)}  &\Gamma_{Z}(Y, \cD(\cS^{*-p}))/\Gamma_{Z\cap Y}(Y, \cD(\cS_Y^{*-p}))      \\
\mapd{\theta'} &     &\mapdr{\theta'} \\
D(\chS^*(X, Y)^\natural)  &\mapr{\cap u}  &\ph{\,.}D(\chS_{Z}^{*-p}(X, Y))\,.
 \end{array}
$$ 
\smallskip 

3.\, One has also a commutative square 
$$
\begin{array}{ccc}
  {\Gamma_{Z}(X, \cD(\cS_Y^*))}/{\Gamma_{Z\cap Y}(Y, \cD(\cS_Y^*))}    &\mapl{}    &  {\Gamma(Z, \cD(\cS_Z^*))}/{\Gamma(Z\cap Y, \cD(\cS_{Z\cap Y}^*))}      \\
\mapd{}&    & \mapd{}     \\
D(\chS_Z^*(X, Y))&\mapl{}  &D(\chS^* (Z, Z\cap Y)) 
\end{array}
$$
with all maps quasi-isomorphism; this is induced from the commutative square in 
(\comparisonprep), item 3.
The second row is the dual of the quasi-isomorphism in item 3 of (\chSrelative).
\smallskip

\newcommand{\sst}{\scriptstyle}
4.\, 
The digram below contains diagram $(\star)_{X, Y}$ we already have, and is the
 relative analogue of the diagram $(\star\star)_{X}$ in (\supportedSheafsingularcomparison):
$$\!\!\!\!\!\!(\star\star)_{X, Y}\qquad\quad
\begin{array}{ccccccc}
 &&\frac{\scriptstyle \Gamma(X, \cD(\cS^*))}{\Gamma(Y, \cD(\cS_Y^*)) } &\mapr{\scriptstyle{\cap \theta(u)}}&\frac{\sst\Gamma_Z(X, \cD(\cS^{*-p}))}{\Gamma_{Z\cap Y}(Y, \cD(\cS_Y^{*-p})) }&\gets &\frac{\sst \Gamma(Z, \cD(\cS_Z^{*-p}))}{{\Gamma(Z\cap Y, \cD(\cS_{Z\cap Y}^*))} }\\
&\raise1ex\hbox{$\scriptstyle{\theta'}$}\!\!\swarrow & \mapd{\sst{\theta'} }&     &\mapdr{\theta'} &&\mapdr{\theta'}\\
{\sst D(\check{S}^*(X, Y))}&\leftarrow&{\sst D(\check{S}^*(X, Y)^\natural)   }&\mapr{\cap u}  &{\sst D(\check{S}_Z^{*-p}(X, Y)) }
&\gets &{\sst D(\check{S}^{*-p}(Z, Z\cap Y) )  }\\
\mapu{\xi}&&&     &&&\mapur{\xi} \\
{\sst \ti{S}_*(X, Y) } &&&  && &{\sst \ph{\,,}\ti{S}_{*-p}(Z, Z\cap Y)\,.}
 \end{array}
$$  
The preceding argument confirms it commutes.  
\bigskip

\noindent (\ThmlfandBMpaircapprod) {\bf Theorem.}\quad{\it 
We have the analogue of Theorem (\supportedSheafsingularcomparison) for a pair of 
spaces $(X, Y)$. 
In particular, one has for $u\in S^p(X, X-Z)$, $du=0$,  a commutative diagram on homology
$$\begin{array}{ccc}
 H_m(X, Y) &\mapr{\cap [\theta(u)]}  &H_{m-p}(Z, Z\cap Y)  \\
 \mapd{} &     &\mapdr{} \\
H^{lf}_m(X, Y)  &\mapr{\cap [u]} &\ph{\,.}H^{lf}_{m-p}(Z, Z\cap Y)\,.
 \end{array}
$$  
with vertical maps canonical isomorphisms.
}\bigskip

We also generalize (\CScapproduct). 
Suppose $X$ is a simplicial complex satisfying the condition mentioned at the beginning of 
\S \Sectsimplicial.  Let $Z$ be a subcomplex of $Z$, $N$ be its closed star;
then $N\cap Y$ is the closed star of $Z\cap Y$ in $Y$. 
Assume that the inclusion $(Z, Z\cap Y)\subset (N, N\cap Y)$ 
is a proper deformation retract. 

For a cocycle $u\in C^p(X, N^c)$ recall that we have the map 
$\cap u: \tiC_*(X)\to \tiC_{*-p}(N)$. 
For the restriction $u\in C^p(Y, N^c\cap Y)$ we have the map 
$\cap u: \tiC_*(Y)\to \tiC_{*-p}(N\cap Y)$.  Therefore there is an induced map
$$\cap u: \tiC_*(X, Y)\to \tiC_{*-p}(N, N\cap Y)\,.$$
It induces a map on homology
$$\cap [u]: H^\simp_m(X, Y)\to H^\simp_{m-p}(N, N\cap Y)\,.$$
Composing with the isomorphism 
$H^\simp_{m-p}(Z, Z\cap Y)\to H^\simp_{m-p}(N, N\cap Y)$ we obtain a map
$$\cap [u]: H^\simp_m(X, Y)\to H^\simp_{m-p}(Z, Z\cap Y)\,.$$

The proof of the next result is parallel to that of (\CScapproduct).
\bigskip 

\noindent (\ThmlfandBMpairsimp) {\bf Theorem.}\quad{\it 
For an element $u\in S^p(X, X-Z)$ with $du=0$, the diagram in the derived categoy
$$\begin{array}{ccc}
\tiS_*(X, Y)    &\mapr{\cap u}  &\tiS_{*-p}(N, N\cap Y)  \\
\mapu{\rho} &     &\mapur{\rho} \\
  \tiC_*(X, Y)  &\mapr{\cap u }  &\tiC_{*-p}(N, N\cap Y)
 \end{array}
$$  
induces a commutative diagram on homology 
$$\begin{array}{ccccc}
H^{lf}_m(X, Y)    &\mapr{\cap [u]} &H^{lf}_{m-p}(N, N\cap Y )  &\overset{\sim}\gets  &H^{lf}_{m-p}(Z, Z\cap Y )  \\
\mapu{\rho} &     &\mapur{\rho} & &\mapur{\rho} \\
H^\simp_m(X, Y)  &\mapr{\cap [u] }  &H^\simp_{m-p}(N, N\cap Y) &\overset{\sim}\gets     &\ph{\,}H^\simp_{m-p}(Z, Z\cap Y)\,.
 \end{array}
$$  
}
\bigskip

{\bf Localization isomorphisms for Borel-Moore homologies}\smallskip

In the rest of this section, we study the canonical isomorphism, called the {\it localization 
isomorphism\,}
$H_*(X, Y)\cong H_*(X-Y)$ (for both the Borel-Moore and locally finite homology)
and its compatibility with cap product. 
For the sheaf theoretic Borel-Moore homology, the existence of such an isomorphism
is a special case of [Br-V, (5.10)]; below we will give a direct proof.

For the locally finite homology, the existence of the localization isomorphism 
requires a few steps. 

Unlike the complex $\Gamma(X, \cD(\cS^*))$, the complex $\tiS_*(X)$ is not 
contravariantly functorial for inclusions of open sets, but there is such functoriality 
in the derived category. 
For an open set $U$ of $X$, we produce maps of complexes 
$$\tiS_*(X)\to \tiS_*(U)_X \gets \tiS_*(U)\,,$$
the latter map being an isomorphism. 
By definition, 
$$\tiS_*(U)_X = \varprojlim_{K\in \cpt(U)} S_*(X, X-K), $$
the inverse limit of the complexes $ S_*(X, X-K)$ for $K\in \cpt(U)$. 
The natural maps $S_*(U, U-K)\to S_*(X, X-K)$ induce a map 
$\tiS_*(U) \to \tiS_*(U)_X$, which will be shown to be a quasi-isomorphism. 
Also there is obviously a natural map $\tiS_*(X)\to\tiS_*(U)_X$. 
If $V$ is an open set contained in $U$, there is a natural map 
$\tiS_*(U)_X\to \tiS_*(V)_X$. 
\bigskip

\noindent (\tiSUSUX) {\bf Proposition.}\quad{\it
The map $\tiS_*(U) \to \tiS_*(U)_X$ is a quasi-isomorphism.
}\smallskip 

\begin{proof} Each map $S_*(U, U-K)\to S_*(X, X-K)$ is a quasi-isomorphism by the excision theorem. 
There is an increasing sequence of compact sets of $U$
$$K_1\subset K_2\subset \cdots $$
such that $K_n\subset Int(K_{n+1})$ and $\cup K_n=U$. 
(For example see [Br 2], Theorem 12.11 and its proof.)
Then $\{C_*^{(n)}=S_*(X, X-K_n)\}_{n\ge 1}$ is an inverse system of complexes with surjective transition maps.
One thus has an exact sequence
$$0\to \varprojlim_n\!{}^1\,H_{m+1}(C_*^{(n)}) \to H_{m}(\varprojlim_n C_*^{(n)})\to \varprojlim_n H_{m}(C_*^{(n)})\to 0\,.$$
The same holds for the inverse system of complexes  $\{S_*(U, U-K_n)\}_{n\ge 1}$, so there results another short 
exact sequence.
The groups $H_{m}(C_*^{(n)})$ are isomorphic for the two inverse systems.
Comparing the two short exact sequences we obtain the assertion.
\end{proof}
\bigskip 

Thus the maps $\tiS_*(X)\to\tiS_*(U)_X\gets \tiS_*(U)$ give a map 
$\tiS_*(X)\to \tiS_*(U)$ in the derived category, and induces on homology a map $H^{lf}_m(X)\to H^{lf}_m(U)$.  

The complex $\tiS_*(X)$ is contravariantly functorial for inclusions of open sets. 
Let $V\subset U\subset X$ be open sets. 
The quasi-isomorphism $\tiS_*(V)\to \tiS_*(V)_X$ factors as the composition of 
quasi-isomorphisms $\tiS_*(V)\to \tiS_*(V)_U \to \tiS_*(V)_X$, 
the latter the inverse limit of the maps $S_*(U, U-K)\to S_*(X, X-K)$ for $K\in \cpt(V)$.
One has a commutative diagram 
\[
  \xymatrix@R=10pt{
\tiS_*(X) \ar[r]  &\,\tiS_*(U)_X \ar[r] & \,\tiS_*(V)_X   \\
 &\tiS_*(U)\ar[u] \ar[r]     &\,\tiS_*(V)_U\ar[u]    \\
 &   &  \tiS_*(V) \ar[u]
  }   
\]
from which the contravariant functoriality follows.
\bigskip

\noindent (\localisom) {\bf Localization isomorphisms}
\smallskip

1.\,  The {\it localization map\,} for singular homology is defined as follows. 
Let $Y$ be closed, locally contractible subspace of $X$, and $U=X-Y$. 
Recall by definition $\tiS_*(X, Y)= \tiS_*(X)/ \tiS_*(Y)$.  The composition of the maps
$\tiS_*(Y)\to \tiS_*(X)\to \tiS_*(U)_X$ is zero, which follows from the map 
$S_*(Y)\to S_*(X, X-K)$ being zero for $K\in \cpt(U)$. 
So there is an induced map of complexes 
$\tiS_*(X, Y)\to \tiS_*(U)_X$; this induces on homology a map $H^{lf}_m(X, Y)\to H^{lf}_m(U)$;
we will show below that this is an isomorphism, and call it the {\it localization isomorphism} for the singular homology.
\smallskip

2.\, Likewise, the composition $\chS^*(U)_X\to \chS^*(X) \to \chS^*(Y)$ is zero, since the composition of maps 
$S^*(X, X-K)\to S^*(X)\to S^*(Y)$ is zero.
So we get a commutative diagram of complexes with exact rows
$$
\begin{array}{ccccccc}
0\to &\chS^*(X, Y)  &\mapr{}    &\chS^*(X)    &\mapr{} &\chS^*(Y) &\to 0\\
& \mapu{}&    & \Vert   & &\mapu{}  &\\
0\to &\chS^*(U)_X  &\mapr{}    &\chS^*(X)   &\mapr{} &\chS_Y^*(X) &\ph{\,.}\to 0\,.
\end{array} 
\eqno{(a)}
$$
The vertical arrows are quasi-isomorphisms.
The dual of the quasi-isomorphism $\chS^*(U)_X\to \chS^*(X, Y)$ is also a quasi-isomorphism.
\sskp

3.\,
We give a variant of the map $\xi: \tiS_*(X)\to D(\chS^*(X))$ as follows. 
 For $K\in \cpt(U)$, we have a map $\xi: S_*(X, X-K)\to D(S^*(X, X-K))$; passing to the inverse 
 limit we get a map $\xi: \tiS_*(U)_X\to D(\chS^*(U)_X)$. 
There is a commutative diagram of complexes 
  $$
\begin{array}{ccccc}
D\chS^*(X)    &\mapr{}    & D (\chS^*(U)_X)  &\gets &  D (\chS^*(U)_X)\\
\mapu{\xi}&    & \mapu{\xi}   &&\mapur{\xi}  \\
\tiS_*(X)&\mapr{}   &\tiS_*(U)_X  &\gets &\tiS_*(U)
\end{array}
$$
where the map $D (\chS^*(U)_X)\to D (\chS^*(U)_X) $ is a quasi-isomorphism. 
The map $\xi$ on the right is a quasi-isomorphism by \S 3. 
Therefore $\xi: \tiS_*(U)_X\to D (\chS^*(U)_X)$ is also a quasi-isomorphism. 
Deduced from the left square is another commutative diagram of complexes
  $$
\begin{array}{ccc}
D\chS^*(X, Y)    &\mapr{}    & D (\chS^*(U)_X)   \\
\mapu{}&    & \mapu{}     \\
\tiS_*(X, Y)&\mapr{}   &\ph{\,.}\tiS_*(U)_X  \,.
\end{array}
$$
The map $\xi$ on left is a quasi-isomorphism by (\chSrelative), item 6, so the map $\tiS_*(X, Y)\to \tiS_*(U)_X$ is also a 
quasi-isomorphism.
\sskp

4.\, 
Let $\al: \Gamma_c(X, \cS^*)\to \Gamma_c(X, \cS^*_Y)$ be the canonical surjection. 
We have a commutative diagram of complexes with exact rows, with vertical arrows quasi-isomorphisms.
$$
\begin{array}{ccccccc}
0\to &\Ker (\al)   &\mapr{}    &\Gamma_c(X, \cS^*)     &\mapr{\al} &\Gamma_c(X, \cS^*_Y) &\to 0\\
& \mapu{}&    & \Vert  & &\mapu{}  &\\
0\to &\Gam_c(U, \cS^*)  &\mapr{}    &\Gamma_c(X, \cS^*)   &\mapr{} &\Gamma_c(X, \cS^*|_Y) &\to 0\,.
\end{array}
\eqno{(b)}
$$
Taking the dual we have a commutative diagram, 
with horizontal maps isomorphisms and vertical maps quasi-isomorphisms, 
$$
\begin{array}{ccc}
 \Gamma(X, \cD(\cS^*))/ \Gamma(Y, \cD(\cS^*_Y))    &\isoto    & D(\Ker\al)     \\
\mapd{}&    & \mapd{}     \\
\Gamma(X, \cD(\cS^*))/D\Gam_c(Y, \cS^*|_Y)&\isoto  &\Gamma(U, \cD(\cS^*))
\end{array}
\eqno{(b')}
$$
in particular a quasi-isomorphism $ \Gamma(X, \cD(\cS^*))/ \Gamma(Y, \cD(\cS^*_Y))
\to \Gamma(U, \cD(\cS^*))$.  
This induces an isomorphism on homology
$H_m(X, Y)\to H_m(U)$, the {\it localization isomorphism} for the Borel-Moore homology.
\sskp

5.\,
Both (a) and (b) are diagrams of the same type consisting of six terms.
There is a map of diagrams $\theta: (a)\to (b)$, given as follows. 
We have obtained before three quasi-isomorphisms (denoted $\theta$)\, $\chS^*(X)\to \Gamma_c(X, \cS^*)$, 
$\chS^*(Y)\to \Gamma_c(X, \cS^*_Y)$, and $\chS_Y^*(X)\to\Gamma_c(X, \cS^*|_Y)$;
they give maps from four terms of (a) to the corresponding four terms of (b). 
There are unique maps of complexes $\theta: \chS^*(X, Y)\to \Ker (\al)$ and 
$\theta: \chS^*(U)_X\to \Gam_c(U, \cS^*) $ which, along the other maps, form
a map of diagrams from $(a)\to (b)$.  This map $\theta$
 is termwise a quasi-isomorphism.

If we take the dual of the map $\theta: (a)\to (b)$, we get a commutative diagram,
$$
\begin{array}{ccc}
D\chS^*(X)/ D\chS^*(Y)   &\isoto   & D\chS^*(X, Y)     \\
\mapd{}&    & \mapd{}     \\
D\chS^*(X)/ D\chS_Y^*(X)&\isoto  & D\chS^*(U)_X
\end{array}
\eqno{(a')}
$$
with horizontal maps isomorphisms and vertical maps quasi-isomorphisms, 
another commutative diagram 
$$
\begin{array}{ccc}
 \Gamma(X, \cD(\cS^*))/ \Gamma(Y, \cD(\cS^*_Y))    &\isoto    & D(\Ker\al)     \\
\mapd{}&    & \mapd{}     \\
\Gamma(X, \cD(\cS^*))/D\Gam_c(Y, \cS^*|_Y)&\isoto  &D(\chS^*(U)_X) 
\end{array}
\eqno{(b')}
$$
also with horizontal maps isomorphisms and vertical maps quasi-isomorphisms, and a map 
diagrams $\theta': (b')\to (a')$
which is termwise a quasi-isomorphism.
In particular, we have a commutative diagram of quasi-isomorphisms
$$\begin{array}{ccc}
 \Gamma(X, \cD(\cS^*))/ \Gamma(Y, \cD(\cS^*_Y)) &\rightarrow{} &\Gamma(U, \cD(\cS^*)) \\
  \mapd{\theta'} &  &\mapdr{\theta'} \\
  D(\chS^*(X, Y)) &\rightarrow{} &\ph{\,.}D(\chS^*(U)_X) \,.
   \end{array}$$

6.\,
We have obtained a commutative diagram 
$$\begin{array}{ccc}
 \Gamma(X, \cD(\cS^*))/ \Gamma(Y, \cD(\cS^*_Y)) &\rightarrow{} &\Gamma(U, \cD(\cS^*)) \\
  \mapd{\theta'} &  &\mapdr{\theta'} \\
  D(\chS^*(X, Y)) &\rightarrow{} &D(\chS^*(U)_X) \\
  \mapu{\xi} & &\mapur{\xi} \\
  \tiS_*(X, Y) &\rightarrow{}&\ph{\,.}\tiS_*(U)_X\,.
   \end{array}$$
Recall that the vertical quasi-isomorphisms give the canonical isomorphisms on homology.
We have proven the following theorem.
\bigskip

\noindent (\ThmlfandBMpairlocal) {\bf Theorem.}\quad
{\it
There are natural isomorphisms $H_m(X, Y)\to H_m(U)$ and 
 $H^{lf}_m(X, Y)\to H^{lf}_m(U)$.  
 These isomorphisms are compatible with the canonical isomorphisms 
 $H_m(X, Y)\cong H^{lf}_m(X, Y)$ and $H_m(U)\cong H^{lf}_m(U)$ 
  from (\ThmlfandBMpair).
  }\bigskip
  
The localization isomorphisms are compatible with cap product:
\bigskip

\noindent (\Thmlocalcapproduct) {\bf Theorem.}\quad
{\it Let $u\in H_Z^p(X)$ and $u|_U\in H_{Z\cap U}(U)$ be its restriction.
The cap product $\cap u$ on $H_m(X, Y)$ and the cap product 
$\cap (u|_U)$ on $H_m(U)$ are compatible via the localization isomorphism, namely the diagram 
$$
\begin{array}{ccc}
H_m(X, Y)    &\isoto  &H_m(U)      \\
\mapd{\cap u}&    & \mapdr{\cap u}     \\
H_{m-p}(Z, Z\cap Y)&\isoto  & H_{m-p}(Z\cap U)
\end{array}
\eqno{(\Thmlocalcapproduct.1)}
$$
commutes. 
The same holds for the groups $H^{lf}_*(X, Y)$ and $H_*^{lf}(U)$.
}
\sskp

\bpf
For $u\in \Gamma_Z(X, \cS^*)$ closed, the square 
$$\begin{array}{ccc}
 \Gamma(X, \cD(\cS^*))/\Gamma(Y, \cD(\cS_Y^*))   &\mapr{}  & \Gamma(U, \cD(\cS^*))     \\
\mapd{\cap u } &     &\mapdr{\cap  u} \\
\Gamma_{Z}(Y, \cD(\cS^{*-p}))/\Gamma_{Z\cap Y}(Y, \cD(\cS_Y^{*-p})) &\mapr{}  
&\Gamma_{Z\cap U}(U, \cD(\cS^*))
 \end{array}
$$ 
clearly commutes, and the first assertion follows. 

Consider next the square 
$$
\begin{array}{ccc}
H^{lf}_m(X, Y)    &\isoto  &H^{lf}_m(U)      \\
\mapd{\cap u}&    & \mapdr{\cap u}     \\
H^{lf}_{m-p}(Z, Z\cap Y)&\isoto  & \ph{\,.}H^{lf}_{m-p}(Z\cap U)\,.
\end{array}
\eqno{(\Thmlocalcapproduct.2)}
$$
There is the canonical isomorphism from each term of square 
(\Thmlocalcapproduct.1) to the corresponding term in (\Thmlocalcapproduct.2), 
so one obtains a cubical diagram. 
The face (\Thmlocalcapproduct.1) commutes; 
two of the other faces commute by 
the compatibility of cap product with canonical isomorphisms, (\supportedSheafsingularcomparison) and 
(\ThmlfandBMpaircapprod); two of the faces commute by 
the compatibility of localization isomorphisms and the canonical isomorphisms, (\ThmlfandBMpairlocal).
From these follows the commutativity of face (\Thmlocalcapproduct.2).
\epf

{\bf Remark.}\, As the proof shows, the commutativity of the squares come from the commutativity 
of the corresponding squares of complexes in the derived category.
\bigskip

{}\bigskip 

{\bf References}
\bigskip

[Bra] Brasselet, J.-P. :
D\'efinition combinatoire des homomorphismes de Poincar\'e.
Ast\'erisque, tome 82-83 (1981) 71-91. 
\smallskip

[Br] Bredon, G. E. :
 Sheaf theory. Second edition. Graduate Texts in Mathematics, 170. Springer-Verlag, New York, 1997.
 \smallskip

[Br 2] Bredon, G. E. : Topology and geometry. Corrected third printing of the 1993 original. 
 Graduate Texts in Mathematics, 139. Springer-Verlag, New York, 1997.
 \smallskip

[BM] Borel, A. and Moore, J. C. :
Homology theory for locally compact spaces. 
Michigan Math. J.  7 (1960) 137-159. 
\smallskip


[Ha] Hatcher, A. : Algebraic topology. Cambridge Univ. Press, Cambridge, 2002. 
\smallskip

[I] Iversen, B. :  Cohomology of sheaves. Universitext. Springer-Verlag, Berlin, 1986.
\smallskip

[Mu] Munkres, J. :  Elements of algebraic topology. Addison-Wesley, Menlo Park, CA, 1984. 

[Sp] Spanier, E. H. : Algebraic topology.  Springer-Verlag, New York-Berlin, 1981.
\smallskip
%
 
\end{document}